\documentclass[12pt]{article} 
\usepackage{mathtools}
\usepackage{amsmath}
\usepackage{hyperref}
\usepackage{yfonts}
\usepackage{amssymb}
\usepackage{mathrsfs}
\usepackage{tikz-cd}
\usepackage{amsthm}

\newtheorem*{theorem*}{Theorem}
\newtheorem{theorem}{Theorem}[section]
\numberwithin{theorem}{subsection}

\newtheorem*{lemma*}{Lemma}
\newtheorem{lemma}[theorem]{Lemma}
\newtheorem{proposition}[theorem]{Proposition}
\newtheorem*{proposition*}{Proposition}

\numberwithin{definition}{subsection}

\author{Joseph DiCapua \and Victor Kolyvagin}

\begin{document}

\date{}

\title{\textbf{Parametrization of Formal Norm Compatible Sequences}}

\maketitle

\begin{abstract}
\noindent
    We give a classification of power series parametrizing Lubin-Tate trace compatible sequences. This proof answers a question posed in the literature by Berger and Fourquaux. Lubin-Tate trace compatible sequences are a generalization of norm compatible sequences, which arise in Iwasawa theory and local class field theory. The result we prove generalizes the interpolation theorem proved by Coleman in the classical norm compatible sequence case. We also, jointly with Victor Kolyvagin, give a method for finding such series explicitly in certain special cases.
\end{abstract}

\tableofcontents

\pagebreak

\section{Introduction}
In his paper "Division Values in Local Fields," Coleman studies a connection between power series and $p$-adic numbers given by the local analytic theory. Coleman notes how Kummer obtained number theoretic applications from different formal operations on power series. One of these applications led to "Iwasawa's explicit descriptions of the Galois structure of various modules connected with local cyclotomic fields." In particular the central result of Coleman's paper allows one to interpolate elements of certain modules appearing in Iwasawa theory.
\\
\\
In this dissertation we take Coleman's result on power series interpolating norm compatible sequences of $p$-adic numbers, and we consider an analogue involving the operation coming from a Lubin-Tate formal group law. We study the sequences defined in "Iwasawa theory and $F$-analytic Lubin-Tate $(\varphi, \Gamma)$-modules," and we give a classification of the formal Lubin-Tate trace compatible sequences parametrized by power series. In particular this gives an answer to the question posed in remark 3.4.7 in the above paper. 
\\
\\
This study is parallel to the study of classical norm compatible sequences, and there are some similarities and differences. One significant difference is that there are certain form Lubin-Tate trace compatible sequences that cannot be interpolated, which we explain in section 3.2. On the other hand it is still possible to interpolate certain sequences and the set of such interpolated sequences seems to be large as shown in sections 3.1 and 3.3. Composing with the maps found in \cite{BergerFourquaux}, the series produced in this dissertation also interpolate corestriction compatible sequences in the first cohomology of a certain Galois representation. The series constructed here can also be viewed as a generalization of Coleman series when the base field is $\mathbb{Q}_p$ and the formal group law is the multiplicative formal group law.
\\
\\
We begin by fixing some notation that will be used for the dissertation. Let $K$ be a fixed finite extension of $\mathbb{Q}_p$. $\mathcal{O}_K$ will always denote the ring of integers in $K$. We also fix a Lubin-Tate formal group law $F$ over $K$ associated to some choice of uniformizer $\pi$ of $\mathcal{O}_K$. $q$ will always denote the size of the residue field of $K$. We will often write $x \oplus y = x\oplus_F y$ for $F(x,y)$ whenever $x,y$ have positive valuation. For every $a \in \mathcal{O}_K$ there is an associated power series $[a](x)$ which is an endomorphism of $F$. We let $\textfrak{F}_n$ denote the kernel of the endomorphism $[\pi^{n+1}](x)$. We also let $H$ be a fixed, complete unramified extension of $K$, and we write $\varphi$ for the Frobenius element of the Galois group of $H$ over $K$. We consider the tower of field extensions defined by

$$H_n = H(\textfrak{F}_n)$$

\noindent
We also fix a sequence $u = (u_n)$ for the remainder of this paper satisfying $u_0$ is a nonzero element of $\textfrak{F}_0$ and $[\pi](u_{n+1}) = u_n$ for each $n$. These conditions guarantee each $u_n$ generates $\textfrak{F}_n$ as an $\mathcal{O}_K$-module.
\\
\\
We say a sequence $(a_n)$ with each $a_n\in H_n$ is norm compatible if it satisfies 

$$N_{n+1,n}(a_{n+1}) = a_n$$

\noindent
where $N_{m,n}$ denotes the norm map from $H_m$ to $H_n$. Coleman was able to parametrize all norm compatible sequences with power series defined over $\mathcal{O}_H$.
One of the central results of "Division Values in Local Fields" is the following:

\begin{theorem*}
    Let $\alpha = (\alpha_n)$ be a sequence with each $\alpha_n \in \mathcal{O}_{H_n}$ satisfying $N_{H_{n+1}/H_n}(\alpha_{n+1}) = \alpha_n$. Then there exists a unique power series $f_\alpha(x)$ in $\mathcal{O}_H[[x]]$ satisfying

    $$(\varphi^{-n}f_\alpha)(u_n) = \alpha_n$$

    \noindent
    for all $n \ge 0$.
\end{theorem*}

\noindent
The aim of this dissertation is to generalize this idea to Lubin-Tate trace compatible sequences. We have the correct definition of a Lubin-Tate trace compatible sequence from \cite{BergerFourquaux}: let $S$ be the set of all sequences $(x_n)$ where each $x_n$ lives in the maximal ideal of $K_n = K(\textfrak{F}_n)$, and the $x_n$ satisfy the recursive relationship

$$\text{Tr}^{LT}_{K_{n+1}/K_n}(x_{n+1}) = [q/\pi](x_n)$$

\noindent
for all $n \ge 0$. In the above, the operation $\text{Tr}^{LT}_{K_{n+1}/K_n}$ is defined exactly like usual trace from $K_{n+1}$ to $K_n$, except we replace addition with the operation $\oplus_F$. For all $x$ in the maximal ideal of $K_{n+1}$, we have

$$\text{Tr}^{LT}_{K_{n+1}/K_n}(x) = \sigma_1(x) \oplus \sigma_2(x) \oplus \ldots \oplus \sigma_q(x)$$

\noindent
where the set $\{\sigma_i\}$ is the set of all automorphisms in $G(K_{n+1}/K_n)$. The main result of this dissertation is the classification of all series $f(x) \in \mathcal{O}_K[[x]]$ such that the sequence $(f(u_n)) \in S$.
\\
\\
We briefly describe the methods used to construct such series: in section 2.1 we construct a map from the $\mathcal{O}_K$-module of all series $f \in \mathcal{O}_K[[x]]$ for which $f(u_n) \in S$ to the kernel of Coleman's trace operator. We find that the kernel of this map is the endomorphism ring of our formal group law, and the image is a submodule of the kernel of Coleman's trace operator. Note that this map, which will be labeled as $\log_F(\phi(f))$ in this dissertation, also appears in \cite{Coleman1} as a constant multiple of the map $\Theta_{\textfrak{F}}(f)$.
\\
\\
We first establish that the set of all interpolated sequences in $S$, that is the set of all $(x_n) \in S$ such that there exists $f \in \mathcal{O}_K[[x]]$ satisfying $f(u_n) = x_n$, is an $\mathcal{O}_K$-module. We denote this module by $\mathscr{A}$. Then for all $a \in \mathcal{O}_K$ and all $f \in \mathscr{A}$ we have $a$ acts on $f$ by the composition $[a](f(x))$. We also have addition in the $\mathcal{O}_K$-module $\mathscr{A}$ is given by the formal group law operation of $F$, so that $f_1 \oplus_F f_2 = F(f_1(x),f_2(x))$ gives the sum of two series $f_1,f_2 \in \mathscr{A}$.
\\
\\
Coleman's trace operator is defined to be the unique function $\mathscr{L}$ on power series satisfying 

$$\mathscr{L}(f)([\pi](x)) = \sum_{z\in\textfrak{F}_0} f(x\oplus z)$$

\noindent
We label the intersection of the kernel of Coleman's trace operator with $\pi\mathcal{O}_K[[x]]$ by $\mathscr{C}$. We take $\mathscr{C}'$ to be the $\mathcal{O}_K$-submodule of $\mathscr{C}$ consisting of all series $g \in \mathscr{C}$ with $g'(0) = 0$. We define the function $\phi$ on $\mathscr{A}$ by sending $f\in\mathscr{A}$ to the series

$$\phi(f) =  [\pi](f(x)) \ominus f([\pi](x))$$

\noindent
and one can check that $\log_F(\phi(f)) \in \mathscr{C}$. Note also that $\phi$ factors through the quotient $\mathscr{A}/\text{End}(F)$ since $\phi$ is $\mathcal{O}_K$-linear and $\text{End}(F) = \text{ker}(\phi)$.
\\
\\
We show that the image $\log_F\circ \phi$ on $\mathscr{A}$ is exactly $\mathscr{C}'$, and it follows that $\log_F\circ\phi$ gives an isomorphism between $\mathscr{A}/\text{End}(F)$ and $\mathscr{C}'$. Note that in the cyclotomic case ($F = G_m$ and $[\pi] = [p] = (1+x)^p-1$), the map $p^{-1}\log_F \circ \phi$ and the exact sequence obtained using this map appear in "Local Units Modulo Circular Units" by Coleman \cite{Coleman3}.
\\
\\
In order to obtain this isomorphism we check that the map $\phi$ sends an arbitrary series $f \in \mathscr{A}$ to some series $\phi(f) = h(x) \in \pi\mathcal{O}_K[[x]]$ satisfying

$$ \text{Tr}^{\text{LT}}_{K_{n+1}/K_n}(h(u_{n+1})) = 0$$

\noindent
for all $n \ge 0$. We label the $\mathcal{O}_K$-module of all such series $h$ by $\mathscr{D}$. Note that this definition of $\mathscr{D}$ is equivalent to stipulating $h$ satisfies the identity

$$ h(x) \oplus h(x\oplus z_1) \oplus \ldots \oplus h(x\oplus z_{q-1}) = 0$$

\noindent
where the $z_i$ are the $0$-th level torsion points of $F$. We show that $\log_F : \mathscr{D} \xrightarrow{\sim} \mathscr{C}$ is an isomorphism of $\mathcal{O}_K$-modules.
\\
\\
We are able to find the image of $\phi$ in $\mathscr{D}$ first by showing that if $h(x) \equiv 0 \mod \deg 2$ and $\pi \mid h$ then one can solve the equation $\phi(f) = h$. We show any such solution $f$ must live in $\mathscr{A}$. We then show that any $h \in \phi(\mathscr{A})$ must have a trivial linear term. We also construct additional series in $\mathscr{A}$ to show that $\log_F(\phi(\mathscr{A}))$ is exactly $\mathscr{C}'$.
\\
\\
We prove the following:

\begin{theorem*}
The exact sequence 

$$ 0 \rightarrow \text{\normalfont End}(F) \rightarrow \mathscr{A} \rightarrow \mathscr{C}' \rightarrow 0$$

\noindent
splits where the map $\text{\normalfont End}(F) \rightarrow \mathscr{A}$ is given by inclusion and the map $\mathscr{A} \rightarrow \mathscr{C}'$ is given by the composition of $\log_F \circ \phi$.
\end{theorem*}

\noindent
There is a surjective $\mathcal{O}_K$-linear map from $\mathcal{O}_K[[x]]$ to $\mathscr{C}$ constructed in \cite{Coleman2} which we use to determine all series in $\mathscr{C}'$ in section 2.5. This will complete the classification of series in $\mathcal{O}_K[[x]]$ parametrizing Lubin-Tate trace compatible sequence.
\\
\\
For an alternate proof that $\log_F\circ \phi : \mathscr{A} \rightarrow \mathscr{C}'$ is a surjection see Theorem 22 of \cite{Coleman1}. Since $\mathscr{C}' \subseteq A$ and the preimage of $\mathscr{C}'$ under the map $\Theta_{\textfrak{F}}$ is exactly $\mathscr{A}$ we get another proof of the surjection.
\\
\\
Let $M$ denote the kernel of Coleman's trace operator in $\mathcal{O}_K[[x]]$. In section 2.2 we give a description of $M$. Using a map defined in lemma 16 of \cite{Coleman2} we are able to construct a surjection of $\mathcal{O}_K[[x]]$ onto $M$. This surjection allows us to give a countable set of power series which generate $M$.
\\
\\
In section 3.1 we use the $\frac{q}{\pi}$-eigenspace of Coleman's trace operator in order to give another method of constructing series which interpolate sequences in $S$. We use this method to prove the following proposition, which shows the submodule of interpolated sequences in $S$ is large in some sense:

\begin{proposition*}
    Assume that $\pi^3 \mid q$. If $z$ is an arbitrary element of the maximal ideal of $K_n$ then there exists $\l \ge 0$ and $x \in S$ with $x$ interpolated such that $x_n = [\pi^l](z)$.
\end{proposition*}

\noindent
If we denote the submodule of interpolated sequences of $S$ by $S_\text{int}$ then the above proposition says that the map $K\otimes_{\mathcal{O}_K} S_{\text{int}} \rightarrow K_n$ defined by sending $\lambda \otimes x$ to $\lambda \log_F(x_n)$ is surjective.
\\
\\
In section 3.2 we show that the more general Lubin-Tate trace compatible sequences setting is different from the classical case with the multiplicative group by showing that there are certain sequences which cannot be interpolated when $|\frac{q}{\pi}|$ is small enough. In particular we show this happens when $\pi^3 \mid q$. We prove the following lemma:

\begin{lemma*}
    Suppose $(\alpha_i) \in S$ is interpolated, so there exists some power series $f(x) \in \mathcal{O}_K[[x]]$ such that $f(u_i) = \alpha_i$ for each $i$. Then assuming $f$ is not the zero series we get that $\lim_{i\rightarrow \infty}|\alpha_i|$ exists and is nonzero.
\end{lemma*}

\noindent
When $\pi^3 \mid q$ we construct nontrivial sequences $\alpha \in S$ satisfying $\lim_{i\rightarrow \infty} |\alpha_i| = 0$ which cannot be interpolated by the above lemma. Note that it is already known that not all sequences in $S$ can be interpolated. See remark 3.4.7 of \cite{BergerFourquaux}. 
\\
\\
In section 3.3 we construct an injection from the $\mathbb{Z}_p$-module of norm compatible sequences of principal units over the tower $\{K_n\}$ into the kernel of Coleman's trace operator when $\pi$ is any uniformizer of $K$ such that $\pi^n \neq q$ for all integer exponents $n$. Again this shows that the submodule of interpolated sequences in $S$ is large in some sense.
\\
\\
In section 4, which is joint work with Victor Kolyvagin, we show that certain Coleman series can be obtained explicitly from isomorphisms of different formal group laws. Let $f(x) \in \mathcal{O}_K[[x]]$ be a power series such that $f(x) \equiv x^q \mod \pi$ and such that $f(x) \equiv \pi x \mod x^2$. Let $R(\pi, q)$ denote the set of all $g \in \mathcal{O}_K[x]$ such that $g(x)$ is monic of degree $q$, $g(x) \equiv x^q \mod \pi$, and $g(x) \equiv \pi x \mod x^2$. Then there are formal group laws $F_f$ and $F_g$ associated to $f$ and $g$ respectively for any choice of $g \in R(\pi, q)$. We show that the collection of isomorphisms $i_{f,g} : F_f \rightarrow F_g$ such that $i_{f,g}(x) \equiv x \mod x^2$ provide a supply of "explicit" norm compatible systems in the tower of fields $\{K_n\}$.
\\
\\
We then consider the set of norm compatible sequences of principal units over the tower $\{K_n\}$ as a $\mathbb{Z}_p[\Delta]$-module where $\Delta$ is the cyclic subgroup of order $p-1$ in $G(K_0/K)$. If $\psi : \Delta \rightarrow \mu_{p-1}$ is a homomorphism we define $e_\psi \in \mathbb{Z}_p[\Delta]$ in order to obtain a decomposition of the module of principal units into $e_\psi$ eigenspaces. We are able to give a new proof that norm compatible sequences are interpolated for certain special cases, namely when the norm compatible sequence is contained in an $e_\psi$-eigenspace where $\psi$ is not the trivial character.
\\
\\
Note that by combining the explicit series of section 4 with the injection from section 3.3 we can obtain a method for generating explicit series which interpolate sequences in $S$ as follows:
\\
\\
Let $r(x)$ be any explicit series interpolating a norm compatible sequence. For example one can take $r(x) = i_{f,g}(x)$ to be an isomorphism of formal group laws. Then by the injection constructed in section 3.3 we have

$$\log([p^r]([q]_{G_m}(r(x))\ominus_{G_m}r([\pi]_F(x)))) = p^rq\log(r(x)) - p^r\log(r([\pi_F](x)))$$

\noindent
where $r$ is large enough so that

$$[p^r](\pi\mathcal{O}_K[[x]]) \subseteq p \mathcal{O}_K[[x]]$$

\noindent
is an explicit series in the kernel of Coleman's trace operator. If we denote the above series by $s(x)$, we can multiply $s$ by an appropriate power of $[\pi](x)$ if necessary to guarantee $s \in \mathscr{C}'$. We then apply the inverse of the map $\phi$ defined in section 2.1 to the explicit series $\exp_F(s(x))$ to obtain an explicit series which interpolates some sequence in $S$. 
\\
\\
\textbf{Acknowledgments:} I would like to thank my advisor Victor Kolyvagin for his guidance throughout my time at the Graduate Center. I am grateful for all of the mathematics I have learned from him. He picked an excellent problem for me to work on, and I am grateful that he encouraged me to persist and find a solution. The discussions we had proved to be invaluable in finding a solution to the problem and exploring related areas of research. I am also thankful for the time he spent reading the manuscript and providing detailed feedback.
\\
\\
Moreover, I would like to thank the members of my committee, Kenneth Kramer and Vladimir Shpilrain, for taking the time to read the manuscript and providing valuable feedback.
\\
\\
I am grateful for my family, who were all very supportive throughout the course of my studies.

\pagebreak

\section{Proofs}

\subsection{Determing the $\mathcal{O}_K$-module of all interpolated sequences in $S$}

Let $\mathscr{A}$ denote the set of all series $f$ living in $\mathcal{O}_K[[x]]$ with $|f(0)| < 1$ and satisfying the relation

$$\sum^{\text{LT}}_{z\in \textfrak{F}_0}f(x\oplus z) = [\frac{q}{\pi}](f([\pi](x)))$$

\noindent
Here and for the rest of this dissertation $\Sigma^{\text{LT}}$ will denote summation taken with the formal group operation as addition. $q$ will always denote the size of the residue field $\mathcal{O}_K/\pi\mathcal{O}_K$. We must define what we mean by the sum of two series with respect to addition from the formal group law.
\\
\\
Suppose both $f$ and $g$ live in $\mathcal{O}_E[[x]]$ where $E$ is some finite extension of $K$. Suppose also that $|f(0)| < 1$ and $|g(0)| < 1$. We prove that the composition of functions $F(f(x),g(x))$ for $|x| < 1$ is given by a unique power series $h(x)$ contained in $\mathcal{O}_E[[x]]$ and satisfying $|h(0)| < 1$. We will always take $f(x) \oplus g(x)$ to mean the unique power series $h(x)$ satisfying $h(x) = F(f(x),g(x))$ as functions on the disc $|x| < 1$.
\\
\\
Lemma 5.2 implies that the composition of functions $F(f(x),g(x))$ for $|x| < 1$ is given by a unique power series in $\mathcal{O}_E[[x]]$ namely the coefficientwise limit of

$$ \sum_{i+j \le N} a_{i,j}f(x)^ig(x)^j$$

\noindent
exists and is this series. Here we have

$$F(x,y) = \sum_{i,j}a_{i,j}x^iy^j$$

\noindent
Note that this proof guarantees expressions of the form

$$f(x \oplus z_1) \oplus f(x \oplus z_2)$$

\noindent
are defined where $z_i \in \textfrak{F}_0$ and $|f(0)|<1$. This is because whenever $z \in \textfrak{F}_0$ we have $f(x\oplus z) \in \mathcal{O}_{K_0}[[x]]$ and satisfies $|f(0 \oplus z)| = |f(z)| < 1$.
\\
\\
Then the equation at the beginning of this section is equivalent to $f$ interpolating some element of $S$, since both sides of the equation live in $\mathcal{O}_K[[x]]$ and agree at all torsion points iff $f$ interpolates some sequence in $S$. 
\\
\\
Then $\mathcal{O}_K$ acts on $\mathscr{A}$ in the following way: if $a \in \mathcal{O}_K$ and $f \in \mathscr{A}$ then $a$ acting on $f$ gives the series $[a](f(x)) \in \mathscr{A}$. If $f_1$ and $f_2$ are two series in $\mathscr{A}$ their sum is given by $f_1(x) \oplus f_2(x) \in \mathscr{A}$. It is routine to check that the action of $\mathcal{O}_K$ is compatible with the addition and they make $\mathscr{A}$ an $\mathcal{O}_K$-module.
\\
\\
Next we check that the endomorphism ring of $F$, $\text{End}(F)$, is contained in $\mathscr{A}$ as a submodule provided $q > 2$. Take arbitrary $a \in \mathcal{O}_K$, then

$$\sum^{\text{LT}}_{z\in \textfrak{F}_0}[a](x\oplus z) = [q]([a](x)) \oplus \sum^{\text{LT}}_{z\in \textfrak{F}_0}[a](z)$$

\noindent
One can check that

$$\sum^{\text{LT}}_{z\in\textfrak{F}_0}[a](z) = 0$$

\noindent
If $\pi \mid a$ this follows when because $[a](z) = 0$ whenever $z\in\textfrak{F}_0$. If $a$ is a unit the above follows when $q > 2$ because 

$$\sum^{\text{LT}}_{z\in\textfrak{F}_0} z = 0$$

\noindent
When there are more than $2$ elements in the residue field of $\mathcal{O}_K$ we can find some $\lambda$ a unit in $\mathcal{O}_K$ with $\lambda - 1$ also being a unit. It follows for such $\lambda$ that

$$[\lambda](\sum^{\text{LT}}_{z\in\textfrak{F}_0} z) = \sum^{\text{LT}}_{z\in\textfrak{F}_0} [\lambda](z) = \sum^{\text{LT}}_{z\in\textfrak{F}_0} z$$

\noindent
so that 

$$[\lambda-1](\sum^{\text{LT}}_{z\in\textfrak{F}_0} z) = 0$$

\noindent
which is only possible if the sum is zero. We conclude for all $a \in \mathcal{O}_K$ that 

$$\sum^{\text{LT}}_{z\in \textfrak{F}_0}[a](x\oplus z) = [qa](x)$$

\noindent
when $q > 2$. Furthermore if $[a]$ is an endomorphism then we get 

$$[\frac{q}{\pi}]([a]([\pi](x))) = [qa](x)$$

\noindent
The above proves $\text{End}(F) \subseteq \mathscr{A}$. Since $\text{End}(F)$ is also an $\mathcal{O}_K$-module with respect to the same formal group law addition as $\mathscr{A}$, and because it is closed under the action of $\mathcal{O}_K$, we have $\text{End}(F)$ is an $\mathcal{O}_K$ submodule of $\mathscr{A}$.
\\
\\
Recall from the introduction that $\mathscr{C}$ denotes the set of all $g \in \pi\mathcal{O}_K[[x]]$ such that $\mathscr{L}(g) = 0$. We check that $\mathscr{C}$ is also an $\mathcal{O}_K$-module, and we construct maps to show that the quotient $\mathscr{A}/\text{End}(F)$ is isomorphic to the submodule of $\mathscr{C}$ consisting of all $g \in \mathscr{C}$ such that $g'(0) = 0$.
\\
\\
For the $\mathcal{O}_K$-module structure on $\mathscr{C}$ we take addition to be addition of power series, and $\mathcal{O}_K$ acts by multiplication, so $a\cdot g = ag$ is scalar multiplication. It is clear that $\mathscr{C}$ is closed under addition since $\mathscr{L}$ is a linear function. One can also see that $\pi \mid ag$ for $a\in \mathcal{O}_K$ and $g\in \mathscr{C}$, and we also have $\mathscr{L}(ag) = 0$ since $\mathscr{L}$ is linear.
\\
\\
Next let $\mathscr{D}$ denote the set of all $h \in \pi\mathcal{O}_K[[x]]$ satisfying the functional equation

$$\sum^{\text{LT}}_{z\in\textfrak{F}_0}h(x\oplus z) = 0$$

\noindent
Then $\mathscr{D}$ is also an $\mathcal{O}_K$-module where we take addition to be given by $F$ and action of $a\in\mathcal{O}_K$ to be composition with the series $[a]$, so that $a\cdot h = [a](h(x))$ for all $a \in \mathcal{O}_K$ and all $h \in \mathscr{D}$. We define an $\mathcal{O}_K$-module map $\phi : \mathscr{A} \rightarrow \mathscr{D}$ by the following:

$$\phi ( f ) = [\pi](f(x))\ominus(f([\pi](x)))$$

\noindent
It is elementary to check that $\phi(f)$ satisfies

$$ \sum^{\text{LT}}_{z\in\textfrak{F}_0} \phi(f)(x \oplus z) = 0$$

\noindent
given that $f \in \mathscr{A}$. It will then follow that $\phi(f) \in \mathscr{D}$ if we show $\pi \mid \phi(f)$. To do this we consider the expression

$$[\pi](f(x)) \ominus f([\pi](x)) \mod \pi$$

\noindent
Note that the first term satisfies $[\pi](f(x)) \equiv f(x)^q \mod \pi$. The second term satisfies $f([\pi](x)) \equiv f(x^q) \mod \pi$. It follows that $\phi(f) \equiv f(x)^q \ominus f(x^q) \mod \pi$. However since $q$ is exactly the size of $\mathcal{O}_K/\pi\mathcal{O}_K$ we must have $f(x)^q \equiv f(x^q) \mod \pi$. From this it follows that $\pi \mid \phi(f)$.
\\
\\
It is elementary to check that $\phi$ is a map of $\mathcal{O}_K$-modules. We also show that the kernel of $\phi$ is exactly the $\mathcal{O}_K$-submodule of $\mathscr{A}$ given by the endomorphisms of $F$.
\\
\\
Suppose we have some $f \in \mathscr{A}$ such that

$$\phi(f) = [\pi](f(x))\ominus f([\pi](x)) = 0$$

\noindent
We show this is only possible if $f \in \text{End}(F)$. Note that $[\pi](f(0)) = f(0)$. If $|f(0)| < 1$ this is only possible if $f(0) = 0$, and since $f\in \mathscr{A}$ we know $|f(0)| < 1$, implying $f(0) = 0$. Now either $f(u_i) = 0$ for all $i$ (in which case $f=0$), or there exists a torsion point of smallest index $i_0$ such that $f(u_{i_0}) \neq 0$. However note that $[\pi](f(u_{i_0})) = f([\pi](u_{i_0})) = f(u_{i_0-1}) = 0$. 
\\
\\
The above is only possible if $f(u_{i_0}) \in \textfrak{F}_0$. The relation $[\pi](f(x)) = f([\pi](x))$ then implies recursively that $f(u_{i_0+n}) \in \textfrak{F}_n\backslash\textfrak{F}_{n-1}$. It also implies $[\pi](f(u_{i_0+n})) = f(u_{i_0+n-1})$. Then the sequence $b_n = f(u_{i_0+n})$ for $n \ge 0$ satsfies $b_n \in \textfrak{F}_n/\textfrak{F}_{n-1}$ and $[\pi](b_n) = b_{n-1}$. Since the Galois group $G_\infty$ acts transitively on such sequences there exists some automorphism $\sigma$ of $K_\infty/K$ such that $\sigma(u_i) = b_i$, so there exists $u_\sigma \in \mathcal{O}_K$ such that $[u_\sigma](u_n) = b_n = f(u_{i_0+n})$. 
\\
\\
It follows from the above equality that $f = [\pi^{i_0}u_\sigma]$ since both series agree on almost all torsion points. Another way of seeing this is to note if two series $f,g \in \mathcal{O}_K[[x]]$ agree on infinitely many values $x_i$ with each $|x_i| < 1$ then $f(x) = g(x)$. For a proof of this see lemma 5.1. This completes the proof that the kernel of $\phi$ is contained in the endomorphisms of $F$. Showing that every endomorphism is in the kernel is elementary, so we conclude that the kernel of $\phi$ is exactly the set of endomorphisms of $F$.
\\
\\
We would now like to study the image $\phi(\mathscr{A})$ in $\mathscr{D}$. We show the submodule of $\mathscr{D}$ consisting of all $h \in \mathscr{D}$ with $h'(0) = 0$ is contained in the image of $\phi$. We prove the following lemma:

\begin{lemma}
Let $g \in \pi\mathcal{O}_K[[x]]$ such that $g'(0) = 0$. Then there exists $f$ in $\mathcal{O}_K[[x]]$ with $|f(0)| < 1$ such that $[\pi](f(x))\ominus(f([\pi](x))) = g(x)$.
\end{lemma}

\noindent
To solve the above we must have $[\pi](f(0))\oplus i_F(f(0)) = g(0)$. This is equivalent to the equality $[\pi -1](f(0)) = g(0)$. Since $\pi-1$ is a unit in $\mathcal{O}_K$ we can take $f(0) = [\frac{1}{\pi-1}](g(0)) \in \pi\mathcal{O}_K$. We now define a sequence of coefficients $(a_n)$, with $a_0 = f(0) = [\frac{1}{\pi-1}](g(0))$, such that if we let $f_N = \sum_{n=0}^N a_nx^n$ then 

\begin{equation}
\label{Cong}
    [\pi](f_N(x))\oplus i_F(f_N([\pi](x))) \equiv g(x) \mod x^{N+1}
\end{equation}

\noindent
Here $i_F = i \in \mathcal{O}_K[[x]]$ is defined to be the unique series satisying $F(x,i_F(x)) = 0$. Once we have the above for all $N$ we can take $f = \sum_{n=0}^\infty a_nx^n$, and $f$ will be a solution to

$$[\pi](f(x))\ominus f([\pi](x)) = g(x)$$

\noindent
Since $b_1 = 0$ where $g(x) = \sum_{n=0}^\infty b_nx^n$ we take $a_1 = 0$ to obtain $f_1$. Now we assume we have solved for all coefficients through $a_N$ and show we can solve for $a_{N+1}$. In what follows we will use derivatives in order to simplify certain expressions. Let $[\pi]'(x)$ denote the first derivative of the series $[\pi](x)$. Let $i'(x)$ denote the first derivative of the series $i(x)$. We need to figure out the coefficient of $x^{N+1}$ in 

\begin{equation}
\label{F-dif}
    [\pi]( f_N(x) + a_{N+1}x^{N+1}) \ominus (f_N([\pi](x)) + a_{N+1}([\pi](x))^{N+1})
\end{equation}

\noindent
We have congruences (mod $x^{N+2})$:

$$[\pi](f_N(x) + a_{N+1}x^{N+1}) \equiv [\pi](f_N(x)) + [\pi]'(f_N(x))a_{N+1}x^{N+1}$$

\noindent
which is equivalent to

$$
[\pi](f_N(x)) + [\pi]'(f_N(0))a_{N+1}x^{N+1} \mod x^{N+2}$$

\noindent
We also have

$$i(f_N([\pi](x)) + a_{N+1}([\pi](x))^{N+1}) \equiv i(f_N([\pi](x))) + i'(f_N([\pi](x)))a_{N+1}([\pi](x))^{N+1}$$

\noindent
which is equivalent to 

$$i(f_N([\pi](x))) + i'(f_N(0))a_{N+1}\pi^{N+1}x^{N+1} \mod x^{N+2}$$

\noindent
Now let $A = [\pi](f_N(x))$, and let $B = i(f_N([\pi](x)))$. Let $\delta = [\pi]'(a_0)a_{N+1}x^{N+1}$, and let $\gamma = i'(a_0)a_{N+1}\pi^{N+1}x^{N+1}$. Then we have:

$$F(A + \delta, B + \gamma) \equiv F(A,B) + \frac{\partial F}{\partial x}(A,B)\delta + \frac{\partial F}{\partial y}(A,B) \gamma $$

\noindent
modulo the ideal generated by $\delta^2$, $\gamma^2$, and $\delta \gamma$. So $(2)$ is equivalent to

$$F(A,B) + \frac{\partial F}{\partial x}([\pi](a_0), i(a_0))[\pi]'(a_0)a_{N+1}x^{N+1} + \frac{\partial F}{\partial y}([\pi](a_0),i(a_0))i'(a_0)a_{N+1}\pi^{N+1}x^{N+1}$$

\noindent
modulo $x^{N+2}$. The congruence $(1)$ will hold for $N=0$ if we put $a_1 = 0$ because $g'(0) = 0$, and if $N \ge 1$, then the coefficient of $a_{N+1}x^{N+1}$ in $(2)$ will be $\pi + \pi^2z$, $z \in \mathcal{O}_K$ because $\frac{\partial F}{\partial x}(0,0) = 1$, $[\pi]'(a_0) = \pi + \pi^2\omega$ with $\omega \in \mathcal{O}_K$. Taking into account $F(A,B) \equiv 0 \mod \pi$, we can determine $a_{N+1}$ as $\frac{C_{N+1}}{(1+\pi z)}$, where $g(x) -F(A,B) = \pi C_{N+1}x^{N+1} \mod x^{N+2}$.
\\
\\
By the above lemma if $h \in \mathscr{D}$ and $h'(0) = 0$ we get there exists an $f\in\mathcal{O}_K[[x]]$ satisfying $[\pi](f)\oplus i_F(f([\pi])) = h$. Since $h$ satisfies the identity

$$ \sum^{\text{LT}}_{z\in\textfrak{F}_0}h(x\oplus z) = 0$$

\noindent
we must have $f \in \mathscr{A}$. To see this we expand $h = \phi(f)$ to get

$$\sum_{z\in\textfrak{F}_0}^{\text{LT}}\phi(f)(x\oplus z) = 0$$

\noindent
so that

$$\sum_{z\in\textfrak{F}_0}^{\text{LT}}[\pi](f(x\oplus z)) \ominus \sum^{\text{LT}}_{z\in\textfrak{F}_0}f([\pi](x\oplus z)) = 0$$

\noindent
the above implies

$$\sum_{z\in\textfrak{F}_0}^{\text{LT}}[\pi](f(x\oplus z)) = [q](f([\pi](x)))$$

\noindent
which is only possible if 

$$\sum_{z\in\textfrak{F}_0}^{\text{LT}}f(x\oplus z) = [q/\pi](f([\pi](x)))$$

\noindent
which is exactly the equation defining $\mathscr{A}$. This last equality follows because the power series $[\pi](x)$ has a formal power series inverse in $K[[x]]$. It follows that if $f_1$ and $f_2$ are two series in $\mathcal{O}_K[[x]]$ with $[\pi](f_1) = [\pi](f_2)$ we must have $f_1 = f_2$ in $\mathcal{O}_K[[x]]$.
\\
\\
Next we use the logarithm and exponential of our formal group law $F$ to show $\mathscr{C}$ and $\mathscr{D}$ are isomorphic as $\mathcal{O}_K$-modules.

\begin{lemma}
The map $\log_F : \mathscr{D} \rightarrow \mathscr{C}$ where $h\in \mathscr{D}$ is sent to the composition $\log_F(h(x))$ is an isomorphism of $\mathcal{O}_K$-modules with inverse given by $\exp_F$.
\end{lemma}

\noindent
Proof: it is well known that $\log_F : \pi\mathcal{O}_K \rightarrow \pi\mathcal{O}_K$ and $\exp_F : \pi \mathcal{O}_K \rightarrow \pi \mathcal{O}_K$ are inverse isomorphisms of $\mathcal{O}_K$-modules. See for example Proposition 7.17 and Proposition 2.4 in \cite{Kolyvagin}, agreeing with the $\mathcal{O}_K$-action follows if we consider Theorem 2 in section 5.1 of \cite{BorevichShafarevich}.The same estimates of divisibility of $\log_F(a)$ and $\exp_F(b)$ depending on divisibility of $a,b$ in the above proof imply that $\log_F : \pi\mathcal{O}_K[[x]] \rightarrow \pi\mathcal{O}_K[[x]]$ and $\exp_F : \pi\mathcal{O}_K[[x]] \rightarrow \pi\mathcal{O}_K[[x]]$ are defined (as coefficientwise limits, see the beginning of section 2.1). The remaining claims follow because they are free for substitutions $x\in\pi\mathcal{O}_K$, and coefficientwise limits agree with composition of functions on $\pi\mathcal{O}_K$, and $\log_F(\mathscr{D}) = \mathscr{C}$.
\\
\\
Remark: for the composition of $\log_F(x)$ with a series $f \in \mathcal{O}_K[[x]]$ to be well defined it suffices that $\pi \mid f(0)$. This follows after taking into account the above mentioned estimates of divisibility of terms of the series $\log_F$ and $\exp_F$.
\\
\\
We then get that $\log_F\circ\phi$ is an $\mathcal{O}_K$-module map from $\mathscr{A}$ to $\mathscr{C}$ which has kernel equal to kernel of $\phi$. It follows immediately that the image $\log_F\circ \phi(\mathscr{A})$ contains $\mathscr{C}'$. This is because $\phi(\mathscr{A})$ contains the submodule of $\mathscr{D}$ consisting of all $h(x)$ with $h'(0) =0$. For arbitrary $g \in \mathscr{C}$, $\exp_F(g)$ is in $\phi(\mathscr{A}) \subseteq \mathscr{D}$ if $g \in \mathscr{C}'$.
\\
\\
We need another lemma to show $\phi(\mathscr{A})$ is exactly the submodule of $\mathscr{D}$ consisting of $h(x) \in \mathscr{D}$ such that $h'(0) = 0$.

\begin{lemma}
Let $f\in \mathcal{O}_K[[x]]$, $|f(0)| < 1$. Then there is a unique endomorphism of $F$, $[\lambda](x)$, such that the linear term of $f(x) \ominus [\lambda](x)$ is trivial.
\end{lemma}

\noindent 
Proof: first we need for such series $f$ the composition $\log_F(f(x))$ is a series in $K[[x]]$ for which evaluation agrees with function composition. For the proof see the remark after lemma 2.1.2.
\\
\\
Next consider $\log_F(f(x))$. $\log_F'(f(0))f'(0) \in \mathcal{O}_K$ since $\log_F'(x) \in \mathcal{O}_K[[x]]$. Then there exists a unique $\lambda \in \mathcal{O}_K$ such that $\log_F(f(x)) - \lambda\log_F(x)$ has trivial linear term. It follows that $f(x) \ominus [\lambda](x)$ has trivial linear term. This is because if we have

$$f(x) \ominus [\lambda](x) = \sum_{n=0}^\infty a_nx^n$$

\noindent
then the linear term of $\log_F(f(x)) - \lambda\log_F(x)$ is given by $\log_F'(a_0)a_1$ and $|a_0| < 1$ so $\log_F'(a_0)$ cannot be zero. It follows that the linear term of $f(x) \ominus [\lambda](x)$ is trivial if and only if the linear term of $\log_F(f(x)) - \lambda\log_F(x)$ is trivial.
\\
\\
From this lemma it follows that if $f \in \mathscr{A}$ there exists a unique $\lambda \in \mathcal{O}_K$ such that the series $f\ominus[\lambda](x)$ has trivial linear term. Then it follows from the following lemma that $\phi(f(x)\ominus[\lambda](x)) \in \mathscr{D}$ also has trivial linear term.

\begin{lemma}
    If $f(x)$ in $\mathcal{O}_K[[x]]$ satisfies $|f(0)| < 1$ and $f'(0) = 0$, then $\phi(f)'(0) = 0$.
\end{lemma}

\noindent
Proof: $f(0) \equiv a_0 \mod x^2$ where $a_0$ is the constant term of $f$ when $f'(0) = 0$. It follows that $[\pi](f) \equiv [\pi](a_0) \mod x^2$. We also have that $f([\pi](x)) \equiv a_0 \mod x^2$. It follows that

$$\phi(f) \equiv [\pi](f(x)) \ominus f([\pi](x)) \equiv [\pi](a_0) \ominus a_0 \mod x^2$$

\noindent
This implies $\phi(f)'(0) = 0$.
\\
\\
Next we note that $[\lambda](x) \in \text{ker}(\phi)$, which implies $\phi(f(x)) = \phi(f(x) \ominus [\lambda](x))$. We conclude that $\phi(\mathscr{A})$ is exactly the submodule of $\mathscr{D}$ consisting of series $h(x) \in \mathscr{D}$ such that $h'(0) = 0$.
\\
\\
It follows that $\log_F\circ \phi(\mathscr{A})$ is exactly $\mathscr{C}'$. Then $\log_F\circ\phi : \mathscr{A} \rightarrow \mathscr{C}'$ factors through an isomorphism of $\mathscr{A}/\text{End}(F)$ with $\mathscr{C}'$. This also just amounts to saying the sequence 

$$ 0 \rightarrow \text{End}(F) \rightarrow \mathscr{A} \rightarrow \mathscr{C}' \rightarrow 0$$

\noindent
where the map $\text{End}(F) \rightarrow \mathscr{A}$ is given by inclusion and the map $\mathscr{A} \rightarrow \mathscr{C}'$ is given by $\log_F\circ \phi$ is exact.
\\
\\
Using lemma 2.1.3 together with the short exact sequence we get the following:

\begin{theorem}
The short exact sequence given by

$$ 0 \rightarrow \text{\normalfont End}(F) \rightarrow \mathscr{A} \rightarrow \mathscr{C}' \rightarrow 0$$

\noindent
splits. Here the map $\text{\normalfont End}(F) \rightarrow \mathscr{A}$ is given by inclusion and the map $\mathscr{A} \rightarrow \mathscr{C}'$ is given by the composition of $\log_F\circ \phi$.
\end{theorem}

\noindent
We define the map $h : \mathscr{A} \rightarrow \text{End}(F)$ by sending $f$ to the unique $[\lambda_f](x) = h(f)$ such that $f\ominus [\lambda_f](x)$ has trivial linear term. We show the map $h$ is a map of $\mathcal{O}_K$-modules. 
\\
\\
We must show $h(f_1\oplus f_2) = [\lambda_{f_1}](x)\oplus [\lambda_{f_2}](x)$, and we must show $h([a](f)) = [a\lambda_f]$. To get the first equality note that $f_1\ominus[\lambda_{f_1}]$ and $f_2\ominus[\lambda_{f_2}]$ both have trivial linear term, so that $F( f_1\ominus[\lambda_{f_1}], f_2\ominus[\lambda_{f_2}] )$ must also have trivial linear term. To see this note that $F(f_1\ominus[\lambda_{f_1}], f_2\ominus[\lambda_{f_2}] )$ must be congruent to a constant mod $x^2$ if $f_1\ominus[\lambda_{f_1}]$ and $f_2\ominus[\lambda_{f_2}]$ both have trivial linear term. Since $[\lambda_{f_1}](x)\oplus [\lambda_{f_2}](x) = [\lambda_{f_1} + \lambda_{f_2}](x)$ it must be the case that $[\lambda_{f_1} + \lambda_{f_2}](x) = h(f_1\oplus f_2)$ by the uniqueness of $\lambda_f$. Then for arbitrary $a\in \mathcal{O}_K$ we have the linear term of $f(x) \ominus [\lambda_f](x)$ is zero, hence the linear term of $[a]( f(x) \oplus [\lambda_f](x) )$ will also be zero. From this it follows that $[\lambda_{[a](f)}] = [a]([\lambda_f](x)) = [a\lambda_f]$. This completes the proof that $h$ is a map of $\mathcal{O}_K$-modules. Then for any $[a](x) \in \text{End}(F) \subset \mathscr{A}$ it is clear that $h([a](x)) = [a](x)$. This completes the proof that the short exact sequence splits.
\\
\\
At this point we note that there is an alternate proof of Theorem 2.1.5 following the arguments in \cite{Coleman1}. As mentioned in the introduction, for another proof that $\log_F\circ \phi : \mathscr{A} \rightarrow \mathscr{C}'$ is a surjection see Theorem 22 of Coleman's paper. First note the map $\Theta_{\textfrak{F}}$ is the same as the map $\pi^{-1} \log_F(\phi(x))$. Since $\pi^{-1}\mathscr{C}' \subseteq A$ provided $q > 2$ and the preimage of $\pi^{-1}\mathscr{C}'$ under the map $\Theta_{\textfrak{F}}$ is exactly $\mathscr{A}$ we get another proof of the surjection.
\\
\\
In summary this section shows that if $f$ is a series in $\mathscr{A}$, then $f$ (up to adding an endomorphism of $F$) corresponds to a series in the kernel of Coleman's trace operator with trivial linear term. In particular if one can find all series in $\mathscr{C}'$ (which is equivalent to finding all series in $\phi(\mathscr{A})$), one can use this isomorphism to find all series $f \in \mathscr{A}$. The next goal of this dissertation is to give a description of $\mathscr{C}'$.
\pagebreak

\subsection{Constructing series in the kernel of Coleman's trace operator}

In this section we will construct series in the kernel of Coleman's trace operator in $\mathcal{O}_K[[x]]$ using the arguments from lemma 16 of \cite{Coleman2}. We will show that the $\mathcal{O}_K$-module, $M =\pi^{-1}\mathscr{C}$, of all such series cannot be finitely generated as an $\mathcal{O}_K$-module. We will also exhibit certain infinite subsets of the kernel which are $\mathcal{O}_K$-linearly independent.
\\
\\
Furthermore in this section we exhibit a countable subset of $M$ which generates all of $M$ by taking possibly infinite sums of series in the subset with coefficients in $\mathcal{O}_K$. We use coefficientwise convergence to show that the infinite sums mentioned above always converge to some series in $\mathcal{O}_K[[x]]$.
\\
\\
The following construction of the series $k(x)$ and $w(x)$ closely follows \cite{Coleman2}. One can find an expression for $\pi/(q-1)$ in $K_0$ as follows. Let $a_i$ be a system of representatives for the residue field of $\mathcal{O}_{K_0}$ consisting of only units and zero. We can further stipulate that all of these representatives live in $\mathcal{O}_K$ since $K_0/K$ is totally ramified. Then one can write:

$$ \pi/(q-1) = \sum_{n=0}^\infty a_{i_n}u_0^n$$

\noindent
since $u_0$ is a uniformizer of $K_0$. Note that we actually have a sum

$$ \pi/(q-1) = \sum_{n= q-1}^\infty a_{i_n}u_0^n$$

\noindent
since $|\pi/(q-1)| = |u_0^{q-1}|$, and any nonzero $a_{i_n}$ with $n < q-1$ would force the right side to have larger absolute value. We define $k(x)$ to be the series

$$k(x) = \sum_{n=q-1}^\infty a_{i_n}x^n$$

\noindent
so that $k(u_0) = \pi/(q-1)$ and $k$ lives in $x^{q-1}\mathcal{O}_K[[x]]$. From this it follows that

$$\mathscr{L}_F(k)(0) = \sum_{z\in \textfrak{F}_0} k(z) = (q-1)\pi/(q-1) = \pi$$

\noindent
(since $k(0) = 0$) which implies $\mathscr{L}_F(k) = \pi w$ for some $w \in \mathcal{O}_K[[x]]$ with $w(0) = 1$. This is also because we know $\pi \mid \mathscr{L}_F(k)$ from lemma 6 of \cite{Coleman1} (see also lemma 5.3).
\\
\\
Now we take arbitrary $g \in \mathcal{O}_K[[x]]$. We consider the expression

$$\mathscr{L}_F(k \frac{g([\pi])}{w([\pi])})$$

\noindent
and we show that this series must equal $\pi g$. To see this note that 

$$\mathscr{L}_F(k \frac{g([\pi])}{w([\pi])})([\pi](x)) = \sum_{z\in\textfrak{F}_0} k(x\oplus z)\frac{g([\pi](x\oplus z)}{w([\pi](x\oplus z)} = \frac{g([\pi](x))}{w([\pi](x))} \sum_{z \in \textfrak{F}_0} k(x\oplus z)$$

\noindent
The expression on the right is just 

$$\frac{g([\pi](x))}{w([\pi](x))} \mathscr{L}_F(k)([\pi](x)) = \frac{g([\pi](x))}{w([\pi](x))} \pi w([\pi](x)) = \pi g([\pi](x))$$

\noindent
which is only possible if 

$$\mathscr{L}_F(k \frac{g([\pi])}{w([\pi])})(x) = \pi g(x)$$

\noindent
This argument shows that the image of $\mathscr{L}_F$ on $\mathcal{O}_K[[x]]$ contains all of $\pi \mathcal{O}_K[[x]]$. By lemma 6 of \cite{Coleman1} (see also lemma 5.3) we have that $\mathscr{L}_F(f) \equiv 0 \mod \pi$ for arbitrary $f \in \mathcal{O}_K[[x]]$. This proves we must have the equality

$$ \mathscr{L}_F(\mathcal{O}_K[[x]]) = \pi \mathcal{O}_K[[x]]$$

\noindent
Note that the above proof is taken from the proof of lemma 16 in \cite{Coleman2}.
\\
\\
Now consider the function $x^n$ for any positive integer $n$ with $n < q-1$. We let $\pi g_n(x) = \mathscr{L}_F(x^n)$ for some choice of $g_n \in \mathcal{O}_K[[x]]$ which is possible since $\mathscr{L}_F(x^n) \equiv 0 \mod \pi$. For this choice of $g_n$ we get 

$$\mathscr{L}_F(k \frac{g_n([\pi])}{w([\pi])}) = \mathscr{L}_F(x^n)$$

\noindent
implying that the series 

$$h_n(x) = x^n - k \frac{g_n([\pi])}{w([\pi])} = x^n - k\frac{\mathscr{L}_F(x^n)([\pi])}{\pi w([\pi])}$$

\noindent
is in the kernel of Coleman's trace operator. Note that each $h_n(x) = x^n + \text{ higher degree terms}$ because $x^{q-1} \mid k(x)$, and because $w(0) = 1$. We know that $w([\pi](x))^{-1} \in \mathcal{O}_K[[x]]$, implying

$$x^{q-1} \mid k(x)\frac{\mathscr{L}_F(x^n)([\pi](x))}{\pi w([\pi](x))}$$

\noindent
so we get that $h_n(x) \equiv x^n \mod x^{n+1}$ because $q-1 > n$.
\\
\\
Now consider the collection of all functions of the form $([\pi](x))^mh_n(x)$ where $m$ ranges over all nonnegative integers, and $n$ is in the range $0 \le n < q-1$. We first check that $([\pi](x))^mh_n(x)$ is also contained in the kernel of $\mathscr{L}_F$. We then use the set of functions $([\pi](x))^mh_n(x)$ to show that the $K$-vector space generated by the kernel of $\mathscr{L}_F$ in $\mathcal{O}_K[[x]]$ is not of finite dimension. We use this fact to show that the kernel of Coleman's trace operator in $\mathcal{O}_K[[x]]$ cannot be finitely generated as an $\mathcal{O}_K$-module.
\\
\\
We must check that $\mathscr{L}_F( ([\pi](x))^mh_n(x) ) = 0$. We have

$$\mathscr{L}_F( ([\pi](x))^mh_n(x) )([\pi](x)) = \sum_{z\in\textfrak{F}_0} ([\pi](x\oplus z))^mh_n(x\oplus z)$$

\noindent
The above is the same as

$$\sum_{z\in\textfrak{F}_0} ([\pi](x\oplus z))^mh_n(x\oplus z) = ([\pi](x))^m \sum_{z\in\textfrak{F}_0} h_n(x\oplus z)$$

\noindent
and this is just $([\pi](x))^m\mathscr{L}_F(h_n(x)) = 0$, which is only possible if 

$$\mathscr{L}_F( ([\pi](x))^mh_n(x) )(x) = 0$$

\noindent
We now consider the vector space $V$ formed by taking all finite $K$-linear combinations of series in $M$. That is 

$$V = \{ \sum_{i=1}^k \lambda_if_i \mid \lambda_i \in K \text{ and each } f_i \in M\}$$

\noindent
We must have the containment of sets $\{ ([\pi](x))^mh_n(x) \}_{m \ge 0, 0\le n < q-1} \subset V$ and we show that the set $\{([\pi](x))^mh_n(x) \}_{m \ge 0, 0 \le n < q-1}$ is linearly independent over $K$.
\\
\\
Suppose there exists a finite linear combination

$$\lambda_1([\pi](x))^{m_1}h_{n_1}(x) + \ldots +\lambda_k([\pi](x))^{m_k}h_{n_k}(x) = 0$$

\noindent
with each $\lambda_i \in K$. Here for distinct $i$ and $j$ we have $(m_i,n_i) \neq (m_j,n_j)$. Wtihout loss of generality assume $|\lambda_1|$ is maximal. Then by multiplying by an appropriate power of $\pi$ if necessary, we can assume each $\lambda_i \in \mathcal{O}_K$ and $\lambda_1$ is a unit. We then consider the above equation mod $\pi$ to get

$$\lambda_1([\pi](x))^{m_1}h_{n_1}(x) + \ldots +\lambda_{k'}([\pi](x))^{m_{k'}}h_{n_{k'}}(x) \equiv 0 \mod \pi$$

\noindent
where we can now assume each remaining $\lambda_i$ is a unit else the term would vanish mod $\pi$. At this point we note that 

$$([\pi](x))^mh_n(x) \equiv x^{qm + n} + \text{ higher degree terms } \mod \pi$$

\noindent
Out of all the pairs $(m_i,n_i)$ appearing in the above equation there must exist a unique $j$ such that $qm_j + n_j$ is minimized. This is because if $m$ is chosen to be minimal, then all of the $n_i$ with $m_i = m$ must be distinct. It follows that for the unique $j$ for which $qm_j+n_j$ is minimal we get 

$$\sum_{i} \lambda_i([\pi](x))^{m_i}h_{n_i}(x) \equiv \lambda_jx^{qm_j + n_j} + \text{ higher degree terms } \mod \pi$$

\noindent
Since $\lambda_j$ is a unit this means the sum is nonzero mod $\pi$ and gives a contradiction if we assume the set $\{ ([\pi](x))^mh_n(x) \}_{m \ge 0, 0 \le n < q-1}$ is linearly dependent.
\\
\\
The above argument shows the $K$-vector space $V$ cannot have finite dimension, and we use this to check that the $\mathcal{O}_K$-module $M$ cannot be finitely generated as an $\mathcal{O}_K$-module.
\\
\\
Suppose that $M$ is finitely generated as an $\mathcal{O}_K$-module. Let $t_1, t_2, \ldots, t_m$ be a generating set for $M$ over $\mathcal{O}_K$. Let 

$$\alpha = \sum_{i=1}^k \lambda_if_i$$

\noindent
be an arbitrary element of $V$, so each $f_i \in M$. Then for each $i$ there exists coefficients $a_{i,n} \in \mathcal{O}_K$ for $1 \le n \le m$ such that

$$f_i = \sum_{n=1}^m a_{i,n}t_n$$

\noindent
These $a_{i,n}$ exist because $f_i \in M$ and $M$ is generated by the $t_n$. Then

$$\alpha = \sum_{i=1}^k (\lambda_i \sum_{n=1}^m a_{i,n}t_n)$$

\noindent
Then the above sum can be written as 

$$\alpha = \sum_{n=1}^m \lambda_n' t_n$$

\noindent
where each $\lambda_n' = \sum_{i=1}^k \lambda_ia_{i,n} \in K$. This implies $\alpha$ is in the $K$ span of the series $t_1,t_2, \ldots, t_m$, so if this is the case then

$$V = \{ \sum_{i=1}^m \lambda_it_i \mid \lambda_i \in K\}$$

\noindent
so that $V$ has finite dimension as a $K$ vector space. However, we have already shown that $V$ cannot be finite dimensional as a $K$-vector space using the series $\{[\pi](x)^mh_n(x)\}_{m\ge 0, 1 \le n \le q-1}$, so our assumption that $M$ is finitely generated as an $\mathcal{O}_K$-module must be false.
\\
\\
In the following paragraphs we show that the short exact sequence

$$0 \rightarrow M \rightarrow \mathcal{O}_K[[x]] \rightarrow \pi\mathcal{O}_K[[x]] \rightarrow 0$$

\noindent
of $\mathcal{O}_K$-modules splits. The map $M \rightarrow \mathcal{O}_K[[x]]$ is given by inclusion, and the map from $\mathcal{O}_K[[x]]$ to $\pi\mathcal{O}_K[[x]]$ is given by $\mathscr{L}_F$.
\\
\\
It will suffice to find a map of $\mathcal{O}_K$-modules 

$$ t : \mathcal{O}_K[[x]] \rightarrow M$$

\noindent
such that $t$ restricted to $M \subset \mathcal{O}_K[[x]]$ is the identity. We claim this map $t$ can be chosen to be 

$$t(g) = g - k \frac{\mathscr{L}_F(g)([\pi](x))}{\pi w([\pi](x))}$$

\noindent
where $k$ and $w$ are the same series defined earlier in this section. It follows that $t$ is linear over $\mathcal{O}_K$ since the map $\mathscr{L}_F$ is linear over $\mathcal{O}_K$. We just need to check $\mathscr{L}_F(t(g)) = 0$ for arbitrary $g \in \mathcal{O}_K[[x]]$ and that $t(g) = g$ if $g \in M$.
\\
\\
To check the first condition note that $\mathscr{L}_F$ is a linear function so that 

$$\mathscr{L}_F(t(g)) = \mathscr{L}_F(g) - \mathscr{L}_F(k \frac{\mathscr{L}_F(g)([\pi](x))}{\pi w([\pi](x))})$$

\noindent
It follows from the argument at the beginning of this section that 

$$\mathscr{L}_F(k \frac{\mathscr{L}_F(g)([\pi](x))}{\pi w([\pi](x))}) = \mathscr{L}_F(g)$$

\noindent
so we must have $\mathscr{L}_F(t(g)) = 0$. To prove that $t$ restricted to $M$ is just the identity it suffices to note that

$$t(g) = g - k \frac{\mathscr{L}_F(g)([\pi](x))}{\pi w([\pi](x))} = g - k\frac{0}{\pi w([\pi](x))} = g$$

\noindent
because $g \in M$ implies $\mathscr{L}_F(g) = 0$. The choice of this map $t$ proves the short exact sequence

$$0 \rightarrow M \rightarrow \mathcal{O}_K[[x]] \rightarrow \pi\mathcal{O}_K[[x]] \rightarrow 0$$

\noindent
splits.
\\
\\
In what follows we give a countable generating set for the kernel of $\mathscr{L}_F$ contained in $\mathcal{O}_K[[x]]$. This set will generate the kernel in that all series in the kernel can be written as (possibly infinite) $\mathcal{O}_K$-linear sums of this set. We use coefficientwise convergence to show that all such linear combinations of elements of the set converge to a series in the kernel.
\\
\\
Let $k$ and $w$ be the same series in $\mathcal{O}_K[[x]]$ defined at the beginning of this section. We extend our definitions of the functions 

$$h_n(x) = x^n - k\frac{\mathscr{L}_F(x^n)([\pi](x))}{\pi w([\pi](x))}$$

\noindent
now to include all integers $n \ge 0$. The same proof used for the cases $n < q-1$ works to show $\mathscr{L}_F(h_n) = 0$. We consider the set of series

$$A = \{ \sum_{n=0}^\infty a_n h_n(x) \mid \text{ each } a_n\in\mathcal{O}_K \}$$

\noindent
We show that each infinite sum of series $\sum_{n=0}^\infty a_nh_n(x) \in A$ defined to be 

$$ \sum_{n=0}^\infty a_nh_n(x) = \lim_{N\rightarrow \infty} \sum_{n=0}^N a_nh_n(x)$$

\noindent
converges coefficientwise to a series in $\mathcal{O}_K[[x]]$.
\\
\\
In order to prove that the above limit converges coefficientwise to some series in $\mathcal{O}_K[[x]]$ it suffices to apply lemma 3.0.2, noting that $h_n(m) \in I_{K_0}^n$ for all $m \in I_{K_0}$ where we take $I_{K_0}$ to be the maximal ideal of $\mathcal{O}_{K_0}$. This concludes the proof that each of the formal sums appearing in the set $A$ converges to some series in $\mathcal{O}_K[[x]]$. We must still show all of these series are in the kernel of Coleman's trace operator. We will use the following lemma:

\begin{lemma}
    Suppose $(f_n(x))$ is a sequence of series in $\mathcal{O}_K[[x]]$ which converges coefficientwise to $f(x)$. Then the sequence of series $(\mathscr{L}_F(f_n))$ converges coefficientwise to $\mathscr{L}(f)$.
\end{lemma}

\noindent
Proof: from the lemma 5.2 it is clear that the sequence of series $(\mathscr{L}_F(f_n)([\pi](x)))$ converges coefficientwise to $\mathscr{L}_F(f)([\pi](x))$ (substituting $m \in I_{K_0}$). This implies by lemma 5.2 that $(\mathscr{L}_F(f_n))$ converges coefficientwise to $\mathscr{L}_F(f)$ because $[\pi] : I \rightarrow I$ is surjective where $I$ is the maximal ideal of $\mathcal{O}_{\overline{K}}$.
\\
\\
From the above lemma we can take any series

$$\sum_{n=0}^\infty a_nh_n(x) \in A$$

\noindent
and it follows that 

$$\mathscr{L}_F(\sum_{n=0}^\infty a_nh_n(x)) = \lim_{N\rightarrow \infty} \mathscr{L}_F(\sum_{n=0}^N a_nh_n(x)) = 0$$

\noindent
since $\mathscr{L}_F(h_n) = 0$ for each $n$.
\\
\\
Conversely, we would like to know that if $f \in \mathcal{O}_K[[x]]$ satisfies $\mathscr{L}_F(f) = 0$ then $f \in A$. This will prove that $A$ is exactly equal to the set of all series $f \in \mathcal{O}_K[[x]]$ such that $\mathscr{L}_F(f) = 0$.
\\
\\
We need to show that if $f$ is an arbitrary element of $\mathcal{O}_K[[x]]$ satisfying $\mathscr{L}_F(f) = 0$ then we have that $f \in A$. We recall the map

$$ t : \mathcal{O}_K[[x]] \rightarrow M$$

\noindent
defined earlier in this section by 

$$t(f) = f - k \frac{\mathscr{L}_F(f)([\pi](x))}{\pi w([\pi](x))}$$

\noindent
We recall the map $t$ is just the identity restricted to $M$, so we have that $t(f) = f$ for $f$ satisfying $\mathscr{L}_F(f) =0 $. We then use this equality to show that if 

$$f(x) = \sum_{n=0}^\infty a_nx^n$$

\noindent
then we also have 

$$ f(x) = \sum_{n=0}^\infty a_n h_n(x)$$

\noindent
when $\mathscr{L}_F(f) = 0$. The equality $t(f) = f$ gives us 

$$ f(x) = \sum_{n=0}^\infty a_nx^n - k \frac{ \mathscr{L}_F( \sum_{n=0}^\infty a_nx^n )([\pi](x))}{\pi w([\pi](x))}$$

\noindent
Since $f(x)$ is the coefficientwise limit of the partial sums $\sum_{n=0}^N a_nx^n$ we can rewrite

$$\mathscr{L}_F(\sum_{n=0}^\infty a_nx^n)([\pi](x))$$

\noindent
in the above equality as

$$\sum_{n=0}^\infty \mathscr{L}_F(a_nx^n)([\pi](x))$$

\noindent
by lemma 2.2.1. Combining the summations on the right side of the equality

$$f(x) = \sum_{n=0}^\infty a_nx^n - \sum_{n=0}^\infty a_nk\frac{\mathscr{L}_F(x^n)([\pi](x))}{\pi w ([\pi](x))}$$

\noindent
gives us exactly the equality

$$f(x) = \sum_{n=0}^\infty a_nh_n(x)$$

\noindent
as desired, so that we see $f \in A$. This allows us to conclude that the set of series $A$ is exactly the set of series in $M$. Stated in other words we get that the kernel of Coleman's trace operator in $\mathcal{O}_K[[x]]$ is equal to the set 

$$ \{ \sum_{n=0}^\infty a_nh_n(x) \mid a_n \in \mathcal{O}_K \}$$

\noindent
Next we assume $q > 2$, and we consider the short exact sequence 

$$0 \rightarrow \mathscr{C}' \rightarrow \mathscr{C} \rightarrow \pi\mathcal{O}_K \rightarrow 0 $$

\noindent
where the map $\mathscr{C}' \rightarrow \mathscr{C}$ is given by inclusion and the map $\mathscr{C} \rightarrow \pi\mathcal{O}_K$ is the projection $g(x) \mapsto g'(0)$. It is clear that the map $\mathscr{C} \rightarrow \pi\mathcal{O}_K$ is a surjection since if we let $g(x) = ah_1(x)$ we have $g'(0) = a$. We will check the given short exact sequence splits.
\\
\\
It suffices to find a map $t : \mathscr{C} \rightarrow \mathscr{C}'$ of $\mathcal{O}_K$-modules such that $t$ restricted to $\mathscr{C}'$ is the identity. We claim the map sending $g(x) \in \mathscr{C}$ to $g(x) - g'(0)h_1(x)$ works as a choice for $t$. First it is clear that

$$t(g) = g(x) - g'(0)h_1(x)$$

\noindent
is linear. We check $t(g) \in \mathscr{C}'$ for arbitrary $g(x) \in \mathscr{C}$. $\mathscr{L}_F(t(g)) = 0$ because $\mathscr{L}_F(g) = \mathscr{L}_F(h_1) = 0$. Also the linear coefficient of $t(g)$ is given by $g'(0) - g'(0)h_1'(0) = 0$ since the linear coefficient of $h_1(x)$ is $1$.
\\
\\
It follows that $t : \mathscr{C} \rightarrow \mathscr{C}'$ is a map of $\mathcal{O}_K$-modules. It then suffices to show that if $g \in \mathscr{C}'$ then $t(g) = g$. This is true for any $g \in \mathscr{C}'$ since we have $g'(0) = 0$ implying $t(g)$ will equal $g$.

\pagebreak

\section{Additional Proofs}

\subsection{The $q/\pi$-eigenspace of Coleman's trace operator}

In this section we use a lemma from "The Arithmetic of Lubin-Tate Division Towers" \cite{Coleman2} to construct elements in the $q/\pi$-eigenspace of Coleman's trace operator. We use these series to give an additional method for constructing series in the module $\mathscr{A}$. We also give a method for finding all series in the $\lambda$-eigenspace of Coleman's trace operator when $\pi$ divides $\lambda$. We do this by constructing an isomorphism between the $\lambda$-eigenspace of Coleman's trace operator and the kernel of Coleman's trace operator, which was already described in section 2.1.
\\
\\
We show that the $\mathcal{O}_K$-module of Lubin-Tate trace compatible sequences $(\alpha_i)$ satisfying $\sup_i|\alpha_i| < 1$ which are interpolated by series $f \in \mathcal{O}_K[[x]]$ is isomorphic to the $\mathcal{O}_K$-module of series in the intersection of the $q/\pi$-eigenspace of Coleman's trace operator with $\pi\mathcal{O}_K[[x]]$. We show the intersection of the $q/\pi$-eigenspace of Coleman's trace operator with $\pi\mathcal{O}_K[[x]]$ is not finitely generated as an $\mathcal{O}_K$-module.
\\
\\
More precisely let $\mathscr{A}'$ be the $\mathcal{O}_K$-module of all series $f(x)$ satisfying $f(x)$ interpolates some sequence $(\alpha_i) \in S$ with $\sup_i|\alpha_i| < 1$. $\mathscr{A}'$ is an $\mathcal{O}_K$-module in the sense that $\lambda \in \mathcal{O}_K$ acts on $f \in \mathscr{A}'$ by $\lambda \cdot f = [\lambda](f(x))$ where $[\lambda](x) \in \text{End}(F)$ is the power series associated to $\lambda$. We also have that addition in $\mathscr{A'}$ is given by $F$, so that for $f_1,f_2 \in \mathscr{A}'$ their sum is given by $F(f_1(x),f_2(x))$.
\\
\\
Let $\mathscr{E}$ denote the set of all series $g \in \pi\mathcal{O}_K[[x]]$ satisfying $\mathscr{L}(g) = \frac{q}{\pi}g$. Then $\mathscr{E}$ is an $\mathcal{O}_K$-module in the sense that $\lambda \in \mathcal{O}_K$ acts on $g\in \mathscr{E}$ by $\lambda \cdot g = \lambda g(x) \in \mathscr{E}$ since $\mathscr{L}$ is linear. Addition in $\mathscr{E}$ is given by addition of power series. Then we have the following:

\begin{theorem}
The map $\log_F : \mathscr{A}' \rightarrow \mathscr{E}$ defined by taking $f(x) \in \mathscr{A}'$ to the composition $g(x) = \log_F(f(x))$ is an isomorphism of $\mathcal{O}_K$-modules. The inverse of this map is given by $\exp_F : \mathscr{E} \rightarrow \mathscr{A}'$ which sends $g(x) \in \mathscr{E}$ to the composition $f(x) = \exp_F(g(x))$.
\end{theorem}

\noindent
We study a property proved in \cite{BergerFourquaux} which shows the set of sequences in $S$ is large in some sense. We show that the set of interpolated sequences in $S$ also has this property. More precisely Berger and Fourquaux showed the following:

\begin{proposition*}
Assume that $K \neq \mathbb{Q}_p$. If $z$ is an arbitrary element of the maximal ideal of $K_n$ then there exists $l \ge 0$ and $x \in S$ such that $x_n = [\pi^l](z)$.
\end{proposition*}

\noindent
We show the following result about the interpolated sequences in $S$ is also true:

\begin{proposition}
Assume that $\pi^3 \mid q$. If $z$ is an arbitrary element of the maximal ideal of $K_n$ then there exists $l \ge 0$ and $x \in S$ with $x$ interpolated such that $x_n = [\pi^l](z)$.
\end{proposition}

\noindent
The additional method of finding series in $\mathscr{A}$ using the $q/\pi$-eigenspace will allow us to prove Proposition 3.1.2.
 \\
 \\
Next we consider the map of $\mathcal{O}_K$-modules $T: \mathcal{O}_K[[x]] \rightarrow \mathcal{O}_K[[x]]$ defined by

$$T(f) = \mathscr{L}(f) - \frac{q}{\pi}f$$

\noindent
It is clear that $f \in \mathscr{E}$ exactly when $T(f) = 0$ and $\pi$ divides $f(x)$. We begin by constructing series in the kernel of $T$ assuming $\pi^3 \mid q$. To construct such nonzero elements in the kernel of $T$ we pick some $f_0 \in \mathcal{O}_K[[x]]$ such that $\mathscr{L}(f_0) = 0$ and such that $f_0 \neq 0 \mod \pi$. One can find such $f_0$ by the description of the kernel of $\mathscr{L}$ given in 2.2. In particular let $h_n$ be the same series from section 2.2. For an arbitrary sequence of coefficients $(a_n)$ with each $a_n \in \mathcal{O}_K$ we let 

$$f(x) = \sum_{n=0}^\infty a_nh_n(x)$$

\noindent
so that $\mathscr{L}(f) = 0$. Then there exists a unique integer exponent $e_f$ such that $\pi^{e_f}f \in \mathcal{O}_K[[x]]$ and $\pi^{e_f}f(x) \neq 0 \mod \pi$. We can then let $f_0(x) = \pi^{e_f}f$. Then 

$$T(f_0) = \mathscr{L}(f_0) - \frac{q}{\pi}f_0 = -\frac{q}{\pi}f_0$$

\noindent
Next we have $\pi^3 \mid q$, so it follows that $-\frac{q}{\pi^2}f_0 \in \pi\mathcal{O}_K[[x]]$. Then by the equality

$$\mathscr{L}(\mathcal{O}_K[[x]]) = \pi \mathcal{O}_K[[x]]$$

\noindent
which can be found in \cite{Coleman2}, we get there exists some $f_1 \in \mathcal{O}_K[[x]]$ such that $\mathscr{L}(f_1) = \frac{q}{\pi^2}f_0$.
\\
\\
It follows that

$$T(f_0 + \pi f_1) = T(f_0) + T(\pi f_1) = -\frac{q}{\pi}f_0 + \mathscr{L}(\pi f_1) - \frac{q}{\pi}\pi f_1 = -qf_1$$
\\
\\
In general suppose we have picked $f_i \in \mathcal{O}_K[[x]]$ for $1 \le i \le N$ such that we know 

$$T(\sum_{i=0}^N\pi^if_i) = -\frac{q}{\pi}\pi^Nf_N$$

\noindent
then we will show we can pick $f_{N+1} \in \mathcal{O}_K[[x]]$ such that 

$$T(\sum_{i=0}^{N+1}\pi^if_i) = T(\sum_{i=0}^N\pi^if_i + \pi^{N+1}f_{N+1}) = -\frac{q}{\pi}\pi^{N+1}f_{N+1}$$

\noindent
is also satisfied. 
\\
\\
Suppose we have chosen $f_i$ up to $i = N$ satisfying the above. Then $\frac{q}{\pi^2}f_N \in \pi\mathcal{O}_K[[x]]$, so we pick $f_{N+1}$ to be any element of $\mathcal{O}_K[[x]]$ such that $\mathscr{L}(f_{N+1}) = \frac{q}{\pi^2}f_N$. It follows from this choice that 

$$\mathscr{L}(\pi^{N+1}f_{N+1}) = \pi^{N+1}\frac{q}{\pi^2}f_N = \frac{q}{\pi}\pi^Nf_N$$

\noindent
and from the above it follows that 

$$T(\sum_{i=0}^N\pi^if_i + \pi^{N+1}f_{N+1}) =$$

$$T(\sum_{i=0}^N\pi^if_i) + T(\pi^{N+1}f_{N+1}) =$$

$$-\frac{q}{\pi}\pi^Nf_N + \mathscr{L}(\pi^{N+1}f_{N+1}) - \frac{q}{\pi}\pi^{N+1}f_{N+1} =$$

$$-\frac{q}{\pi}\pi^{N+1}f_{N+1}$$

\noindent
The sequence of partial sums of the form

$$\sum_{i=0}^N\pi^if_i$$

\noindent
converges coefficientwise, so we can define their limit to be 

$$f = \sum_{i=0}^\infty \pi^if_i$$

\noindent
and $f(x)$ is a well-defined series in $\mathcal{O}_K[[x]]$. Note that for each partial sum we have

$$\sum_{i=0}^N\pi^if_i \equiv f_0 \mod \pi$$

\noindent
so we get that $f \equiv f_0 \mod \pi$ as well. From lemma 2.2.1 it follows that $T(f)$ must equal the coeffcientwise limit of the sequence of series

$$T(\sum_{i=0}^N\pi^if_i) = -\frac{q}{\pi}\pi^Nf_N$$

\noindent
Since this limit is zero we get that $T(f) = 0$, and $f \in \mathcal{O}_K[[x]]$ satisfies 

$$\mathscr{L}(f) = \frac{q}{\pi}f$$

\noindent
Next we will show that the $\mathcal{O}_K$-module of series $f \in \mathcal{O}_K[[x]]$ satisfying

$$\mathscr{L}(f) = \frac{q}{\pi}f$$

\noindent
cannot be finitely generated.
\\
\\
The series $k$ and $w$ were defined in section 2.2. We use the same definitions here. Then the series $h_n(x)$ for $n < q-1$ is still defined to be

$$h_n(x) = x^n - k(x)\frac{\mathscr{L}(x^n)([\pi](x))}{\pi w([\pi](x))}$$

\noindent
We have that $h_n(x) \equiv x^n \mod x^{n+1}$ and $\mathscr{L}(h_n) = 0$. 
\\
\\
We also define $g_{n,m}(x)$ to be 

$$g_{n,m}(x) = [\pi^m](x)h_n(x)$$

\noindent
Then $g_{n,m}(x) \equiv x^{qm+n} + \text{ higher degree terms} \mod \pi$. One can also check that $\mathscr{L}(g_{n,m}) = 0$ as 

$$\mathscr{L}(g_{n,m})([\pi](x)) = \sum_{z \in \textfrak{F}_0} [\pi^m](x\oplus z)h_n(x\oplus z) = [\pi^m](x) \sum_{z\in \textfrak{F}_0} h_n(x\oplus z) = 0$$

\noindent
where we get the last equality because $\mathscr{L}(h_n) =0$. Then by the method described on pages 29-31 we can make a choice of $G_{n,m}(x) \in \mathcal{O}_K[[x]]$ which is associated to $g_{n,m}(x)$ in the sense that $G_{n,m}(x) \equiv g_{n,m}(x) \mod \pi$ and such that

$$\mathscr{L}(G_{n,m})(x) = \frac{q}{\pi}G_{n,m}(x)$$

\noindent
for each pair $n,m$ of nonnegative integers with $n < q-1$.
\\
\\
Now let $M$ be the $\mathcal{O}_K$-module $M = \pi^{-1}\mathscr{E}$, so that $M$ is the module of all series $g(x) \in \mathcal{O}_K[[x]]$ satisfying

$$\mathscr{L}(g)(x) = \frac{q}{\pi}g(x)$$

\noindent
Assume $M$ is finitely generated as an $\mathcal{O}_K[[x]]$-module. Then we can find some nontrivial linear relationship

$$\lambda_{n_1,m_1}G_{n_1,m_1}(x) + \ldots + \lambda_{n_k,m_k}G_{n_k,m_k}(x) = 0$$

\noindent
where $(n_i,m_i) \neq (n_j,m_j)$ for $i \neq j$ and each $\lambda_{n_i,m_i} \in \mathcal{O}_K$. Without loss of generality we may assume $\lambda_{n_1,m_1}$ is a unit after dividing by the correct power of $\pi$ if necessary.
\\
\\
We reduce the above equality mod $\pi$ to get

\begin{equation}
\label{sumGmodpi}
\lambda_{n_1,m_1}G_{n_1,m_1}(x) + \ldots +\lambda_{n_{k'},m_{k'}}G_{n_{k'},m_{k'}} \equiv 0 \mod \pi
\end{equation}

\noindent
where we can now assume each coefficient is a unit, else the term would vanish mod $\pi$. Because each series $G_{n_i,m_i}(x) \equiv g_{n_i,m_i}(x) \equiv x^{qm_i+n_i} + \text{ higher degree terms } \mod \pi$ we get that

$$\lambda_{n_i,m_i}G_{n_i,m_i}(x) \equiv \lambda_{n_i,m_i}x^{qm_i + n_i} + \text{ higher degree terms } \mod \pi$$

\noindent
Because each pair $(n_i,m_i)$ appearing in the sum is unique and each $n_i < q-1$ there must exist a unique index $i$ such that $qm_i + n_i$ is minimal. It follows that for this $i$ we get the left side of $(3)$ is

$$\lambda_{n_i,m_i}x^{qm_i+n_i} + \text{ higher degree terms } \neq 0 \mod \pi$$

\noindent
Therefore we get a contradiction if we assume $M$ is finitely generated as an $\mathcal{O}_K$-module, so we must have that $M$ is not finitely generated.
\\
\\
In the next part of this section we prove Theorem 3.1.1. To prove this theorem we first use that if $f \in \mathscr{A}'$ then $f(u_n) = \alpha_n$ for some sequence $(\alpha_n) \in S$ satisfying $\sup|\alpha_n| < 1$. We use this to show $\pi$ divides $f(x)$ in $\mathcal{O}_K[[x]]$.
\\
\\
Suppose $\pi$ does not divide $f(x)$. Label the coefficients of $f(x)$, so that

$$f(x) = \sum_{n=0}^\infty a_nx^n$$

\noindent
Because $\pi$ does not divide $f$ some of these coefficients must be units. Pick $k_0$ to be the smallest index such that $a_{k_0}$ is a unit. We will show that

$$\lim_{n\rightarrow \infty}|f(u_n)| = \lim_{n\rightarrow \infty} |a_{k_0}u_n^{k_0}| = 1$$

\noindent
First note that if $k$ is any index less than $k_0$ then $\pi \mid a_k$. Then because $\lim_{n\rightarrow\infty}|u_n| = 1$ we will get that there exists some $N_k$ such that whenever $n > N_k$ we have $|a_ku_n^k| < |a_{k_0}u_n^{k_0}|$. Take $N$ to be any integer such that $N > N_k$ for all indices $k < k_0$. We show that for $n > N$ we must have

$$|f(u_n)| = |a_{k_0}u_n^{k_0}|$$

\noindent
For any such $n$ we have

$$|a_ku_n^k| < |a_{k_0}u_n^{k_0}|$$

\noindent
for all $k < k_0$. Also for any index $k > k_0$ we have $|a_k| \le |a_{k_0}|$ and $|u_n^k| < |u_n^{k_0}|$. It follows that

$$|a_ku_n^k| < |a_{k_0}u_n^{k_0}|$$

\noindent
for each $k > k_0$. This is enough to prove that $|f(u_n)|  = |a_{k_0}u_n^{k_0}|$ for sufficiently large indices $n$. This implies that 

$$\sup|f(u_n)| = \sup|\alpha_n| =  1$$

\noindent
which is a contradiction.
\\
\\
We now have that if $f \in \mathscr{A}'$ then $f\in \pi \mathcal{O}_K[[x]]$. By earlier estimates in section 2.1 this allows us to take the composition of series $\log_F(f(x))$. This composition will agree with composition of functions and we also know that $\log_F(f(x)) \in \pi\mathcal{O}_K[[x]]$ by the arguments in 2.1.
\\
\\
If $f \in \mathscr{A}'$ then $f$ satisfies the functional equation

$$\sum_{z\in\textfrak{F}_0}^{\text{LT}}f(x\oplus z) = [\frac{q}{\pi}](f([\pi](x)))$$

\noindent
Both sides of this equation live in $\pi\mathcal{O}_K[[x]]$ and we already saw that $\log_F$ takes addition in the formal group law to ordinary addition of functions in section 2.1. For these reasons if we let $g(x) = \log_F(f(x))$ we get that $g$ will satisfy

$$\sum_{z\in\textfrak{F}_0} g(x\oplus z) = \frac{q}{\pi}g([\pi](x))$$

\noindent
Because the left side of the above equation is $\mathscr{L}(g)([\pi](x))$ the above equation is equivalent to

$$\mathscr{L}(g)(x) = \frac{q}{\pi}g(x)$$

\noindent
Combining all of the above proves that if $f \in \mathscr{A}'$ then the composition $\log_F(f(x))$ lives in $\mathscr{E}$. One can see that composition with $\log_F$ respects the action of $\mathcal{O}_K$ on $\mathscr{A}'$ from the arguments in 2.1. It follows that $\log_F : \mathscr{A}' \rightarrow \mathscr{E}$ is a map of $\mathcal{O}_K$-modules.
\\
\\
One can also check that $\exp_F : \mathscr{E} \rightarrow \mathscr{A}'$ defined by sending $g(x)$ to the composition $\exp_F(g(x))$ is a map of $\mathcal{O}_K$-modules. To show this note that if $g \in \mathscr{E}$ then $\pi$ divides $g(x)$ so that the composition of series $\exp_F(g(x))$ agrees with composition of functions. We also have $\exp_F(g(x)) \in \pi\mathcal{O}_K[[x]]$. Composing series with $\exp_F$ also respects the $\mathcal{O}_K$-module structures of both $\mathscr{E}$ and $\mathscr{A}'$. All of this follows from the arguments presented in section 2.1. 
\\
\\
Now $g(x) \in \mathscr{E}$ exactly if $\pi \mid g(x)$ and 

$$\mathscr{L}(g)(x) = \frac{q}{\pi}g(x)$$

\noindent
This implies we have 

$$\sum_{z\in\textfrak{F}_0}g(x\oplus z) = \frac{q}{\pi}g([\pi](x))$$

\noindent
Applying $\exp_F$ to both sides of the above tells us

$$\sum_{z\in\textfrak{F}_0}^{\text{LT}}f(x\oplus z) = [\frac{q}{\pi}](f([\pi](x)))$$

\noindent
where $f(x)$ equals the composition $f(x) = \exp_F(g(x))$. This tells us that the image of $\exp_F$ on $\mathscr{E}$ is contained in $\mathscr{A}'$. Now since the maps $\log_F$ and $\exp_F$ are inverses by the arguments in section 2.1 we get that both maps are isomorphisms of $\mathcal{O}_K$-modules. This concludes the proof of Theorem 3.1.1.
\\
\\
We now move to the proof of Proposition 3.1.2. For any choice of $l$ we take the logarithm of $[\pi^l](z)$ to get $\pi^l\log_F(z)$. We choose $y_{n-1} \in K_{n-1}$ to be the number such that $y_{n-1}$ satisfies

$$\text{Tr}_{K_n/K_{n-1}}( \log_F(z) ) = \frac{q}{\pi} y_{n-1}$$

\noindent
Then for each $i$ with $0 \le i < n$ we choose the sequence of numbers $y_i \in K_i$, so that the $y_i$ satisfy the relation

$$ \text{Tr}_{K_i/K_{i-1}}(y_i) = \frac{q}{\pi} y_{i-1}$$

\noindent
for $1 \le i < n$. Note that our choice of $y_i$ uniquely determines $y_{i-1}$, so there is only one way to do this. We then pick $l$ large enough so that $\pi^l\log_F(z) \in \pi\mathcal{O}_{K_n}$, and so that $\pi^ly_i \in \pi^{n+1-i}\mathcal{O}_{K_i}$ for each $i$. It is possible to find such an $l$ sufficiently large because we only have finitely many values $y_i$. Then by lemma 9 of \cite{Coleman1}, see also 5.4, we have that there exists a series $f \in \mathcal{O}_K[[x]]$ such that $f(u_i) = \pi^ly_i$ for each $i$ with $i < n$, and such that $f(u_n) = \pi^l\log_F(z)$, and also such that $f(0) = 0$. 
\\
\\
At this point note that 

$$\text{Tr}_{K_n/K_{n-1}}( f(u_n) ) = \text{Tr}_{K_n/K_{n-1}}( \pi^l\log_F(z) ) = \frac{q}{\pi} \pi^l y_{n-1} = \frac{q}{\pi} f(u_{n-1})$$

\noindent
and also that 

$$ \text{Tr}_{K_i/K_{i-1}}(f(u_i)) = \text{Tr}_{K_i/K_{i-1}}( \pi^l y_i) = \frac{q}{\pi} \pi^l y_{i-1} = \frac{q}{\pi} f(u_{i-1})$$

\noindent
for each $i$ with $1 \le i < n$. The above equalities implies the series 

$$\mathscr{L}(f)([\pi](x)) - \frac{q}{\pi} f([\pi](x))$$

\noindent
has zeros at each torsion point $u_i$ for $1 \le i \le n$. We must adjust $f$ so that the above series is also zero at $x = u_0$. Note that

$$\mathscr{L}([\pi](u_0)) - \frac{q}{\pi} f([\pi](u_0)) = 0$$

\noindent
is the same as writing

$$f(0) + \sum_{z\in\textfrak{F}_0}f(z) - \frac{q}{\pi}f(0) = 0$$

\noindent
so it suffices to find a series $F(x) \in \mathcal{O}_K[[x]]$ with $F(u_i) = f(u_i)$ for $0 \le i \le n$ and also such that

$$F(0) = -(1-\frac{q}{\pi})^{-1}\text{Tr}_{K_0/K}(f(u_0))$$

\noindent
Note that $\pi^{n+1}$ divides the above constant because $f(u_0) = \pi^ly_0 \in \pi^{n+1}\mathcal{O}_{K_0}$. It then follows from lemma 5.4 that there exists a series $f_0(x) \in \mathcal{O}_K[[x]]$ such that $f_0(0) = 0$ and such that

$$f_0(u_i) = (1-\frac{q}{\pi})^{-1}\text{Tr}_{K_0/K}(f(u_0))$$

\noindent
for each $i$ with $0 \le i \le n$.
\\
\\
At this point we define the series $F(x) \in \mathcal{O}_K[[x]]$ to be

$$F(x) = f(x) - (1-\frac{q}{\pi})^{-1}\text{Tr}_{K_0/K}(f(u_0)) + f_0(x)$$

\noindent
Then for this choice of $F$ one has $F(u_i) = f(u_i)$ for each $0 \le i \le n$ and

$$F(0) = -(1-\frac{q}{\pi})^{-1}\text{Tr}_{K_0/K}(f(u_0)) = -(1-\frac{q}{\pi})^{-1}\text{Tr}_{K_0/K}(F(u_0))$$

\noindent
Then from the above arguments we see that

$$\mathscr{L}(F)([\pi](x)) - \frac{q}{\pi}F([\pi](x))$$

\noindent
still has zeros at each $u_i$ for $1 \le i \le n$ because $F(u_i) = f(u_i)$ for $0 \le i \le n$. Now because $F$ has the correct constant term we get that $u_0$ is also a zero of the above series. We also get that $0$ is a zero of the above series.
\\
\\
Note that if $u_i$ satisfies the equality

$$\mathscr{L}(F)([\pi](u_i)) = \frac{q}{\pi}F([\pi](u_i))$$

\noindent
then all of the conjugates of $u_i$ also satisfy the equation. This is because

$$\sigma(\mathscr{L}(F)([\pi](u_i))) = \sigma(\frac{q}{\pi}F([\pi](u_i)))$$

\noindent
so that

$$\mathscr{L}(F)([\pi](\sigma(u_i))) = \frac{q}{\pi}F([\pi](\sigma(u_i)))$$

\noindent
for an arbitrary automorphism $\sigma \in G(K_n/K)$.  We conclude that $[\pi^{n+1}](x)$ divides the series $\mathscr{L}(F)([\pi](x)) - \frac{q}{\pi} F([\pi](x))$ in $\mathcal{O}_K[[x]]$.
\\
\\
Define the series $g(x)$ to be 

$$ g(x) = \mathscr{L}(F)([\pi](x)) - \frac{q}{\pi} F([\pi](x)) \in \mathcal{O}_K[[x]]$$

\noindent
We wish to find another series $h(x) \in \mathcal{O}_K[[x]]$ with $[\pi^{n+1}](x)$ dividing $h$ in $\mathcal{O}_K[[x]]$ and satisfying 

$$\mathscr{L}(h)([\pi](x)) - \frac{q}{\pi} h([\pi](x)) = g(x)$$

\noindent
If we can find such a series  $h$ then the series $F-h$ will satisfy $(F-h)(u_n) = \pi^l\log_F(z)$ and will also satisfy

$$\mathscr{L}(F-h)([\pi](x)) - \frac{q}{\pi} (F-h)([\pi](x)) = 0$$

\noindent
so that $(F-h)(x)$ will interpolate some sequence $x_i$ satisfying 

$$\text{Tr}_{K_i/K_{i-1}}(x_i) = \frac{q}{\pi}x_{i-1}$$

\noindent
for all indices $i \ge 1$.
\\
\\
It suffices to find a series $h \in \mathcal{O}_K[[x]]$ with $[\pi^{n+1}](x)$ dividing $h(x)$ such that

$$ \mathscr{L}(h) - \frac{q}{\pi}h = \mathscr{L}(F) - \frac{q}{\pi}F$$

\noindent
Since $[\pi^{n+1}](x)$ divides $g(x)$, we know that $[\pi^n](x)$ divides the right side of the above equation. Also assuming we have the factor $q/\pi$ is divisible by $\pi$, we get that the right side of the above equation is divisible by $\pi$ because we know $\pi \mid \mathscr{L}(F)$ for $F \in \mathcal{O}_K[[x]]$ by lemma 9 of \cite{Coleman1}. It then suffices to prove the following lemma:

\begin{lemma}
Let $j \in \mathcal{O}_K[[x]]$ be an arbitrary series satisfying $\pi \mid j(x)$. Then one can find a series $h(x) \in \mathcal{O}_K[[x]]$ with $[\pi^{n+1}](x)$ dividing $h(x)$ such that 

$$\mathscr{L}(h) - \frac{q}{\pi}h = [\pi^n](x)j(x)$$
\end{lemma}

\noindent
We have the equality $\mathscr{L}(\mathcal{O}_K[[x]]) = \pi\mathcal{O}_K[[x]]$, so we can find a series $h_1(x)$ satisfying

$$ \mathscr{L}(h_1) = j$$

\noindent
We then consider the series $[\pi^{n+1}](x)h_1(x)$. Note that we have

$$ \mathscr{L}( [\pi^{n+1}](x) h_1(x) )([\pi](x)) = \sum_{z\in\textfrak{F}_0} [\pi^{n+1}](x\oplus z) h_1(x\oplus z)$$

\noindent
The right side of the above is just

$$ [\pi^{n+1}](x) \sum_{z \in \textfrak{F}_0} h_1(x\oplus z) = [\pi^{n+1}](x) \mathscr{L}(h_1)([\pi](x))$$

\noindent
The above is only possible if 

$$\mathscr{L}([\pi^{n+1}](x)h_1(x)) = [\pi^n](x) \mathscr{L}(h_1) = [\pi^n](x) j(x)$$

\noindent
If we use $T$ to denote the map $T: \mathcal{O}_K[[x]] \rightarrow \pi\mathcal{O}_K[[x]]$ defined by 

$$T(f) = \mathscr{L}(f) - \frac{q}{\pi}f$$

\noindent
then we get that

$$ T([\pi^{n+1}](x) h_1(x)) = [\pi^n](x)j(x) - \frac{q}{\pi} [\pi^{n+1}](x) h_1(x)$$

\noindent
Now we are still working under the assumption $\pi^3$ divides $q$, so there exists a series $h_2(x) \in \pi\mathcal{O}_K[[x]]$ satisfying 

$$\mathscr{L}(h_2) = \frac{q}{\pi}h_1$$

\noindent
By the same argument as above for the series $h_1$ and $j$ we replace $h_1$ with $h_2$ and we replace $j$ with $\frac{q}{\pi}h_1$ to get that

$$\mathscr{L}([\pi^{n+2}](x)h_2(x)) = \frac{q}{\pi} [\pi^{n+1}](x)h_1(x)$$

\noindent
From these choices of $h_1$ and $h_2$ it follows that

$$T( [\pi^{n+1}](x)h_1(x) + [\pi^{n+2}](x) h_2(x)) = [\pi^n](x)j(x) - \frac{q}{\pi}[\pi^{n+2}](x)h_2(x)$$

\noindent
Now suppose we have found series $h_i(x)$ for $1 \le i \le N$ such that $h_i(x) \in \pi^{i-1}\mathcal{O}_K[[x]]$ and which also satisfy

$$ T( \sum_{i=1}^N [\pi^{n+i}](x)h_i(x) ) = [\pi^n](x)j(x) - \frac{q}{\pi}[\pi^{n+N}](x)h_N(x)$$

\noindent
We check that we can find a series $h_{N+1}(x) \in \pi^N\mathcal{O}_K[[x]]$ such that

$$T( \sum_{i=1}^{N+1} [\pi^{n+i}](x)h_i(x) ) = [\pi^n](x)j(x) - \frac{q}{\pi}[\pi^{n+N+1}](x)h_{N+1}(x)$$

\noindent
$\pi^3$ divides $q$ and $\pi^{N-1} \mid h_N(x)$ implying $\pi^{N+1} \mid \frac{q}{\pi} h_N(x)$. This means there exists a series $h_{N+1}(x) \in \pi^N\mathcal{O}_K[[x]]$ satisfying $\mathscr{L}(h_{N+1}) = \frac{q}{\pi}h_N$. We have

$$\mathscr{L}([\pi^{n+N+1}](x) h_{N+1}(x))([\pi](x)) = \sum_{z\in\textfrak{F}_0} [\pi^{n+N+1}](x\oplus z)h_{N+1}(x\oplus z)$$

\noindent
The right side of the above is 

$$[\pi^{n+N+1}](x) \sum_{z\in\textfrak{F}_0} h_{N+1}(x\oplus z) = [\pi^{n+N+1}](x) \mathscr{L}(h_{N+1})([\pi](x))$$

\noindent
This is only possible if we have

$$\mathscr{L}( [\pi^{n+N+1}](x) h_{N+1}(x) ) = [\pi^{n+N}](x)\mathscr{L}(h_{N+1})(x) = [\pi^{n+N}](x) \frac{q}{\pi}h_N$$

\noindent
It follows from this equality that

$$T( \sum_{i=1}^{N+1} [\pi^{n+i}](x)h_i(x) ) = [\pi^n](x)j(x) - \frac{q}{\pi}[\pi^{n+N+1}](x)h_{N+1}(x)$$

\noindent
At this point we make the choice of $h$ to be the infinite sum of series

$$h(x) = \sum_{i=1}^\infty [\pi^{n+i}](x)h_i(x)$$

\noindent
Note that since each $h_i(x) \in \pi^{i-1}\mathcal{O}_K[[x]]$ the above infinite sum converges coefficientwise to some series in $\mathcal{O}_K[[x]]$. We already saw that $T$ respects coefficientwise limits earlier in this section, and this implies that $T(h)$ is equal to the limit of the expressions $T$ evaluated at the partial sums of the above sum. Precisely we have

$$ T(h) = \lim_{N\rightarrow\infty} T( \sum_{i=1}^N [\pi^{n+i}](x) h_i(x))$$

\noindent
We have already evaluated $T$ at each partial sum, and we showed that

$$T( \sum_{i=1}^N [\pi^{n+i}](x) h_i(x) ) = [\pi^n](x)j(x) - \frac{q}{\pi}[\pi^{n+N+1}](x)h_{N+1}(x)$$

\noindent
It follows that

$$T(h) = \lim_{N\rightarrow \infty} [\pi^n](x)j(x) - \frac{q}{\pi}[\pi^{n+N+1}](x)h_{N+1}(x)$$

\noindent
Since $\pi^N \mid h_{N+1}(x)$ it follows that $\lim_{N\rightarrow\infty}h_{N+1}(x) =0$. It follows that we must have

$$T(h) = [\pi^n](x)j(x)$$

\noindent
This completes the proof of lemma 3.1.3.
\\
\\
At this point we recall the series $F \in \mathcal{O}_K[[x]]$ satisfying $F(u_n) = \pi^l\log_F(z)$ and $F(u_i) = \pi^ly_i$ for $0 \le i < n$. We know that for our construction of the series $F$ we have $[\pi^n](x)$ divides the series 

$$\mathscr{L}(F) - \frac{q}{\pi}F$$

\noindent
We also have that the above series is $0$ mod $\pi$. It therefore follows from the above lemma that one can find a series $h(x) \in \mathcal{O}_K[[x]]$ such that

$$\mathscr{L}(h) - \frac{q}{\pi}h = \mathscr{L}(F) - \frac{q}{\pi}F$$

\noindent
and such that we also have $[\pi^{n+1}](x)$ divides $h(x)$ in $\mathcal{O}_K[[x]]$. We now consider the series $f_1(x) = F(x) - h(x)$. Note that for this choice of $f_1$ we have 

$$\mathscr{L}(f_1) - \frac{q}{\pi}f_1 = \mathscr{L}(F) - \frac{q}{\pi}F - \mathscr{L}(h) + \frac{q}{\pi}h = 0$$

\noindent
We also have that $f_1(u_n) = F(u_n) = \pi^l\log_F(z)$ and $f_1(u_i) = F(u_i) = \pi^ly_i$ for each $i$ with $0 \le i < n$. This is because we have $h(u_j) = 0$ for each index $j$ with $0 \le j \le n$. We also get that $f_1(x)$ satisfies the following identity:

$$\sum_{z\in\textfrak{F}_0} f_1(x\oplus z) = \frac{q}{\pi}f_1([\pi](x))$$

\noindent
This follows from the previous identity

$$\mathscr{L}_F(f_1) - \frac{q}{\pi}f_1 = 0$$

\noindent
and then expanding

$$\mathscr{L}_F(f_1)([\pi](x)) - \frac{q}{\pi}f_1([\pi](x)) = 0$$

\noindent
This identity for $f_1$ implies the sequence $(f_1(u_i))_{i\ge 0}$ satisfies the recursive relation

$$ \text{Tr}_{K_i/K_{i-1}}( f_1(u_i) ) = \frac{q}{\pi} f_1(u_{i-1})$$

\noindent
for all indices $i \ge 1$. We then consider the series $\pi f_1(x)$ which is guaranteed to live in $\pi\mathcal{O}_K[[x]]$. This series satisfies the same relation, namely that

$$ \text{Tr}_{K_i/K_{i-1}}( \pi f_1(u_i) ) = \frac{q}{\pi} \pi f_1(u_{i-1})$$

\noindent
and now we have that all of the values $\pi f_1(u_i)$ for $i \ge 0$ live inside the disc of convergence of the exponential function of the formal group law $F$. Since $\pi f_1(x)$ lives in $\pi\mathcal{O}_K[[x]]$ we get that the composition of series $\exp_F(\pi f_1(x))$ also lives in $\pi\mathcal{O}_K[[x]]$. We also get that the composition of series $\exp_F(\pi f_1(x))$ evaluated at any $x$ with positive valuation agrees with taking $\exp_F$ of the value $\pi f_1(x)$ since all terms of $\pi f_1(x)$ live in the disc of convergence of $\exp_F$ whenever we have $|x| < 1$. For a more detailed explanation of why this is true see section 2.1.
\\
\\
The conclusion of the above paragraph is that the series $\exp_F(\pi f_1(x))$ interpolates some sequence in $S$. We can see this either by referring to Theorem 3.1.1 or by applying $\exp_F$ to both sides of the equality

$$\text{Tr}_{K_i/K_{i-1}}( \pi f_1(u_i) ) = \frac{q}{\pi} \pi f_1(u_{i-1})$$

\noindent
The resulting equality is 

$$ \text{Tr}^{\text{LT}}_{K_i/K_{i-1}} (\exp_F( \pi f_1(u_i))) = [\frac{q}{\pi}](\exp_F (\pi f_1( u_{i-1} )))$$

\noindent
which is exactly the recursive relationship that defines sequences in $S$. We are just left with checking the value of $\exp_F(\pi f_1(u_n) )$ is of the correct form. We have that $\exp_F( \pi f_1(x)) $ evaluated at $u_n$ is equal to $\exp_F$ evaluated at $\pi f_1(u_n)$. Now $\pi f_1(u_n) = \pi^{l+1}\log_F(z)$. Consider the value $\alpha = [\pi^{l+1}](z)$. Then

$$\log_F(\alpha) = \log_F([\pi^{l+1}](z)) = \pi^{l+1}\log_F(z)$$

\noindent
Since $\alpha$ and $\pi^{l+1}\log_F(z)$ both live in a disc on which $\exp_F$ and $\log_F$ are inverse isomorphisms, it follows that the above is only possible if $\exp_F(\pi^{l+1}\log_F(z)) = \alpha = [\pi^{l+1}](z)$. It follows that $\exp_F(\pi f_1(x) )$ evaluated at $u_n$ must equal $[\pi^{l+1}](z)$, and since we have already shown that the series $\exp_F( \pi f_1(x) )$ interpolates some sequence in $S$, we get that $[\pi^{l+1}](z)$ is the $n$-th entry of some interpolated sequence in $S$ as desired. This completes the proof of Proposition 3.1.2. This shows that the set of all interpolated sequences in $S$ is big in the same sense that the set $S$ is big.
\\
\\
For the remainder of this section take $\lambda \in \mathcal{O}_K$ such that $\pi \mid \lambda$. We show that the $\mathcal{O}_K$-module of series $f \in \mathcal{O}_K[[x]]$ satisfying 

$$\mathscr{L}(f) = \lambda f$$

\noindent
is isomorphic to the kernel of $\mathscr{L}$ in $\mathcal{O}_K[[x]]$. Denote the $\mathcal{O}_K$-module of all $f \in \mathcal{O}_K[[x]]$ satisfying $\mathscr{L}(f) = \lambda f$ by $N_\lambda$. We also denote the kernel of $\mathscr{L}$ in $\mathcal{O}_K[[x]]$ by $M$. We define a map of $\mathcal{O}_K$-modules by sending $f \in N_\lambda$ to the series

$$f(x) - \frac{\lambda k(x) f([\pi](x))}{\pi w([\pi](x))} $$

\noindent
Here the series $k(x)$ and $w(x)$ are the same series from the proof of lemma 16 in \cite{Coleman2}. For the definitions of these series see also section 2.2 pages 16-17. We show the above map is an isomorphism of $\mathcal{O}_K$-modules between $N_\lambda$ and $M$.
\\
\\
We refer to the above map as $\rho_\lambda$ defined on $N_\lambda$. We first check that $\rho_\lambda(f)$ is contained in $M$ for arbitrary $f \in N_\lambda$. It suffices to check that $\rho_\lambda(f) \in \mathcal{O}_K[[x]]$ and also that $\mathscr{L}(\rho_\lambda(f)) = 0$.
\\
\\
We get that $\rho_\lambda(f) \in \mathcal{O}_K[[x]]$ because $\pi \mid \lambda$ and $w([\pi](x))$ is a unit in $\mathcal{O}_K[[x]]$.
\\
\\
Next we have

$$\mathscr{L}(\rho_\lambda(f)) = \mathscr{L}(f) - \mathscr{L}( \frac{\lambda k(x) f([\pi](x))}{\pi w([\pi](x))})$$

\noindent
Since $f \in N_\lambda$ we have the above expression equals

$$\mathscr{L}(\rho_\lambda(f)) = \lambda f - \mathscr{L}( \frac{\lambda k(x) f([\pi](x))}{\pi w([\pi](x))})$$

\noindent
Then since $\mathscr{L}$ is linear we get

$$\mathscr{L}(\rho_\lambda(f)) = \lambda f - \lambda \mathscr{L}( \frac{ k(x) f([\pi](x))}{\pi w([\pi](x))})$$

\noindent
By the construction of the series $k$ and $w$ we know that 

$$\mathscr{L}( \frac{k(x) f([\pi](x))}{\pi w([\pi](x))}) = f(x)$$

\noindent
and putting this together with the previous equation gives that

$$\mathscr{L}(\rho_\lambda(f)) = \lambda f - \lambda f = 0$$

\noindent
We conclude that $\rho_\lambda(f) \in M$ for arbitrary $f \in N_\lambda$. One can check that $\rho_\lambda$ is a map of $\mathcal{O}_K$-modules. We would like to show $\rho_\lambda : N_\lambda \rightarrow M$ is injective. We do this by showing that if $f$ is a series in the kernel of $\rho_\lambda$, so that

$$\rho_\lambda(f) = f(x) - \frac{\lambda k(x) f([\pi](x))}{\pi w([\pi](x))} = 0$$

\noindent
then we must have $f(x) = 0$. If $f$ satisfies the above equation then we must have $x^{q-1} \mid f(x)$ because $x^{q-1} \mid k(x)$ implying

$$\frac{\lambda k(x) f([\pi](x))}{\pi w([\pi](x))} \equiv 0 \mod x^{q-1}$$

\noindent
Let $N$ be the largest integer for which $x^N$ divides $f(x)$, which must exist if $f(x) \neq 0$. We get a contradiction if we assume $N$ exists.
\\
\\
If $x^N$ divides $f(x)$ then we have $x^N$ divides $f([\pi](x))$ because $x \mid [\pi](x)$. Then because $x^{q-1} \mid k(x)$ and $w([\pi](x))$ is a unit in $\mathcal{O}_K[[x]]$ it follows that 

$$\frac{\lambda k(x) f([\pi](x))}{\pi w([\pi](x)} \equiv 0 \mod x^{N+q-1}$$

\noindent
The above implies $x^{N+q-1}$ divides $f(x)$, which is a contradiction since $N + q-1 > N$ where we assumed $N$ is the largest integer such that $x^N$ divides $f(x)$. We conclude that the kernel of $\rho_\lambda$ must be trivial, and that $\rho_\lambda : N_\lambda \rightarrow M$ is an injective map of $\mathcal{O}_K$-modules.
\\
\\
We still need to show $\rho_\lambda$ is surjective. Let $h(x) \in M$ be any series in $M$. We construct $g(x) \in N_\lambda$ such that $\rho_\lambda(g) = h$.
\\
\\
In order to find $g$ we will construct a sequence of series $(g_i)$ for which the partial sums

$$\sum_{n=0}^Ng_n(x)$$

\noindent
converge coefficientwise. We will pick the $g_i$ such that the series

$$\rho_\lambda(\sum_{n=0}^Ng_n(x))$$

\noindent
converge coefficientwise to $h$, and this will be enough to show $\rho_\lambda(g) = h$ where 

$$g(x) = \sum_{n=0}^\infty g_n(x)$$

\noindent
Pick $g_1(x) = h(x)$. Then we have

$$\rho_\lambda(g_1) = \rho_\lambda(h) = h(x) - \frac{\lambda k(x) h([\pi](x))}{\pi w([\pi](x))}$$

\noindent
In particular $x^{q-1}$ divides $\rho_\lambda(g_1)-h$. We take 

$$g_2(x) = h(x) - \rho_\lambda(g_1(x)) = \frac{\lambda k(x) h([\pi](x))}{\pi w ([\pi](x))}$$

\noindent
and for this choice of $g_2$ we see that

$$\rho_\lambda(g_2) = \frac{\lambda k(x) h([\pi](x))}{\pi w([\pi](x))} - \frac{\lambda k(x) g_2([\pi](x))}{\pi w([\pi](x))}$$

\noindent

\noindent
Now we have that

$$\rho_\lambda(g_1 + g_2) = h(x) - \frac{\lambda k(x) g_2([\pi](x))}{\pi w([\pi](x))}$$

\noindent
In the above $x^{q-1}\mid k(x)$ and $x^{q-1} \mid g_2([\pi](x))$ implying $x^{2(q-1)}$ divides $\rho_\lambda(g_1 + g_2) - h(x)$.
\\
\\
Now take some integer $N \ge 2$. Suppose for all $1 \le n \le N$ we have picked a sequence $g_n(x) \in \mathcal{O}_K[[x]]$ such that $x^{(n-1)(q-1)}$ divides $g_n$ and such that

$$\rho_\lambda(\sum_{n=0}^N g_n(x)) - h(x)$$

\noindent
is divisible by $x^{N(q-1)}$. We show we can choose $g_{N+1}(x)$ such that $x^{N(q-1)}$ divides $g_{N+1}$ and such that $x^{(N+1)(q-1)}$ divides

$$\rho_\lambda(\sum_{n=0}^{N+1}g_n(x)) - h(x)$$

\noindent
It suffices to let $g_{N+1}(x)$ be

$$g_{N+1}(x) = h(x) - \rho_\lambda(\sum_{n=0}^Ng_n(x))$$

\noindent
For this choice of $g_{N+1}(x)$ note we immediately have $x^{N(q-1)}$ divides $g_{N+1}(x)$. We also have

$$\rho_\lambda(g_{N+1}(x)) = h(x) - \rho_\lambda(\sum_{n=0}^Ng_n(x)) - \frac{\lambda k(x) g_{N+1}([\pi](x))}{\pi w([\pi](x))}$$

\noindent
Because $x^{q-1}$ divides $k(x)$ and $x^{N(q-1)}$ divides $g_{N+1}(x)$ it follows that

$$\rho_\lambda(\sum_{n=0}^{N+1}g_n(x)) = h(x) - \frac{\lambda k(x) g_{N+1}([\pi](x))}{\pi w([\pi](x))}$$

\noindent
and the difference 

$$\rho_\lambda(\sum_{n=0}^{N+1}g_n(x)) - h(x)$$

\noindent
is divisible by $x^{(N+1)(q-1)}$. This completes the proof that we can find $g_{N+1}(x)$ satisfying the necessary conditions.
\\
\\
Now since $x^{(n-1)(q-1)}$ divides $g_n(x)$ we get that the partial sums

$$\sum_{n=0}^Ng_n(x)$$

\noindent
converge coefficientwise to some series $g(x) \in \mathcal{O}_K[[x]]$. By lemma 4.0.2 the expressions

$$\rho_\lambda( \sum_{n=0}^N g_n(x))$$

\noindent
will also converge coefficientwise to

$$\rho_\lambda(g) = g(x) - \frac{\lambda k(x) g([\pi](x))}{\pi w([\pi](x))}$$

\noindent
Because we already showed that $x^{N(q-1)}$ divides

$$\rho_\lambda(\sum_{n=0}^Ng_n(x)) - h(x)$$

\noindent
we get that we must have $\rho_\lambda(g) = h$.
\\
\\
We are left with showing $g(x) \in N_\lambda$. We apply $\mathscr{L}$ to both sides of 

$$\rho_\lambda(g) = g(x) - \frac{\lambda k(x) g([\pi](x))}{\pi w([\pi](x))} = h(x)$$

\noindent
to get 

$$\mathscr{L}(g) - \mathscr{L}(\frac{\lambda k(x) g([\pi](x))}{\pi w([\pi](x))}) = \mathscr{L}(h) = 0$$

\noindent
It follows that

$$\mathscr{L}(g) - \lambda \mathscr{L}( \frac{k(x) g([\pi](x))}{\pi w([\pi](x))}) = 0$$

\noindent
Then because

$$\mathscr{L}(\frac{k(x) g([\pi](x))}{\pi w([\pi](x))}) = g(x)$$

\noindent
by the arguments in lemma 16 of \cite{Coleman2} we have $\mathscr{L}(g) - \lambda g = 0$. We conclude that $g \in N_\lambda$ for every $h \in M$, and we also conclude the map $\rho_\lambda : N_\lambda \rightarrow M$ is an isomorphism of $\mathcal{O}_K$-modules. The above argument shows that for $\pi \mid \lambda$ the $\lambda$-eigenspace of $\mathscr{L}$ is isomorphic to the kernel of $\mathscr{L}$ contained in $\mathcal{O}_K[[x]]$. This allows us to construct series in the $\lambda$-eigenspace of $\mathscr{L}$.

\pagebreak

\subsection{There exist sequences in $S$ which are not interpolated}

We begin this section by showing the following lemma:

\begin{lemma}
    Suppose $(\alpha_i) \in S$ is interpolated, so there exists some power series $f(x) \in \mathcal{O}_K[[x]]$ such that $f(u_i) = \alpha_i$ for each $i$. Then assuming $f$ is not the zero series we get that $\lim_{i\rightarrow \infty}|\alpha_i|$ exists and is nonzero.
\end{lemma}

\noindent
Proof: label the coefficients of $f$ so that

$$f(x) = \sum_{i=0}^\infty a_ix^i$$

\noindent
where each $a_i \in \mathcal{O}_K$. Because each $a_i \in \mathcal{O}_K$ we have there exists at least one index $j$ such that $|a_j| \ge |a_i|$ for each $i$. Take $j_0$ to be the minimal such index $j$, so that $|a_{j_0}| \ge |a_i|$ for each $i$, and so that if $j$ is any other index with $|a_j| \ge |a_i|$ for each $i$ then $j_0 \le j$. We will show that $\lim_{i\rightarrow \infty}|f(u_i)| = |a_{j_0}|$.
\\
\\
We will show the above by showing that

$$|f(u_i)| = |a_{j_0}u_i^{j_0}|$$

\noindent
for sufficiently large values of $i$.
\\
\\
First assume $j < j_0$. We would like to compare the absolute values of the terms $a_ju_i^j$ and $a_{j_0}u_i^{j_0}$. Here it suffices to note that $|u_i| \rightarrow 1$ as $i \rightarrow \infty$. Then because $|a_{j_0}| > |a_j|$ from the definition of $j_0$ we get that

$$|a_{j_0}u_i^{j_0}| > |a_j| \ge |a_ju_i^j|$$

\noindent
for sufficiently large indices $i$. Note here that if $j_0 = 0$ then this case does not apply.
\\
\\
Now we switch to the case $j > j_0$. For such indices and for any value of $i$ we have

$$|a_{j_0}u_i^{j_0}| = |a_{j_0}||u_i^{j_0}| \ge |a_j||u_i^{j_0}| > |a_j||u_i^j| = |a_ju_i^j|$$

\noindent
From the above inequalities we get that

$$|f(u_i)| = |a_{j_0}u_i^{j_0}|$$

\noindent
for sufficiently large indices $i$. It follows from this equality that the following limit exists and we have

$$\lim_{i\rightarrow\infty}|\alpha_i| = \lim_{i\rightarrow \infty} |f(u_i)| = |a_{j_0}|$$

\noindent
because $\lim_{i\rightarrow \infty}|u_i^{j_0}| = 1$. This completes the proof of the lemma.
\\
\\
On the other hand we can also show that if $|q/\pi|$ is small enough then there exist nontrivial sequences $(\alpha_i) \in S$ satisfying $\lim_{i\rightarrow\infty}|\alpha_i| = 0$. By lemma 3.2.1 any such sequence cannot be interpolated. This shows that in general when $|q/\pi|$ is small enough there exist compatible sequences in $S$ not interpolated by power series, so the general case differs from the classical situation with the multiplicative formal group law. 
\\
\\
We now move to the proof that if $|q/\pi|$ is small enough then there exist nontrivial sequences $(\alpha_i) \in S$ satisfying $\lim_{i\rightarrow \infty}|\alpha_i| = 0$. Here we take $|q/\pi|$ small enough to mean $\pi^3 \mid q$.
\\
\\
We will find for each $\epsilon > 0$ a nontrivial sequence $\alpha_\epsilon = (\alpha_{\epsilon,n})$ such that $\sup_n|\alpha_{\epsilon,n}| < \epsilon$, $\alpha_\epsilon \in S$, and $\alpha_\epsilon$ is not interpolated because $\lim_{n\rightarrow \infty}|\alpha_{\epsilon,n}| = 0$. 
\\
\\
Let $r > 0$ be some radius for which $\exp_F$ and $\log_F$ are inverse isomorphisms on the disc $|x| < r$. The existence of such an $r$ follows from the discussion in the proof of lemma 2.1.2 in section 2.1.
\\
\\
We can also stipulate $r$ be small enough so that $|\exp_F(x)| < q^{-(q-1)^{-1}} = |u_0|$ whenever $|x| < r$. This in particular guarantees that $|\log_F(\exp_F(x))| = |\exp_F(x)| = |x|$. For this equality see V of \cite{Coleman1}. We begin by constructing the sequence $(\log_F(\alpha_{\epsilon,n}))$. 
\\
\\
Take any initial value $a_0 \in K_0$ with $|a_0| < \min{r,\epsilon}$. Then we construct a sequence $(a_n)$ with $a_n\in K_n$ such that $\text{Tr}_{K_{n+1}/K_n}(a_{n+1}) = \frac{q}{\pi} a_n$ for all $n\ge 0$ satisfying $|a_n| \le |\frac{q}{\pi^2}|^n|a_0|$ for each $n > 0$. 
\\
\\
Suppose we have found a sequence satisfying the above conditions up to the $N-1$st entry, it is sufficient to find $a_N$ satisfying the conditions and $|a_N| \le |\frac{q}{\pi^2}a_{N-1}|$. We know that $\text{Tr}_{K_N/K_{N-1}}(\mathcal{O}_{K_N}) = \pi \mathcal{O}_{K_{N-1}}$. For this see Proposition 3.4.4 in \cite{BergerFourquaux}. Then pick some $z\in \mathcal{O}_{K_N}$ with $\text{Tr}_{K_N/K_{N-1}}(z) = \pi$. Then we take $a_N = z\frac{q}{\pi^2}a_{N-1}$. For this choice of $a_N$ we have 

$$\text{Tr}_{K_N/K_{N-1}}(a_N) = \text{Tr}_{K_N/K_{N-1}}(z\frac{q}{\pi^2}a_{N-1}) = \frac{q}{\pi^2}a_{N-1}\text{Tr}_{K_N/K_{N-1}}(z) = \frac{q}{\pi}a_{N-1}$$

\noindent
It is then also true that 

$$|a_N| = |z| |\frac{q}{\pi^2}| |a_{N-1}| \le |\frac{q}{\pi^2}||a_{N-1}| \le |\frac{q}{\pi^2}| |\frac{q}{\pi^2}|^{N-1} | a_0| = |\frac{q}{\pi^2}|^N|a_0|$$

\noindent
Choosing the entries in this way proves we have a sequence $(a_n)$ with $a_n\in K_n$ satisfying $\text{Tr}_{K_{n+1}/K_n}(a_{n+1}) = \frac{q}{\pi}a_n$ with $|a_n| \le |\frac{q}{\pi^2}|^n|a_0|$. Also $a_0$ was chosen such that $|a_0| < \epsilon$, so we get that $\sup_n |a_n| < \epsilon$.
\\
\\
At this point we show $\alpha_n$ defined to be $\exp_F(a_n)$ is in $S$ and cannot be interpolated. We have that the sequence $\alpha_n$ defined in this way satisfies $|\alpha_n| = |a_n|$, so that $\sup_n|\alpha_n| < \epsilon$. To show $(\alpha_n) \in S$ it suffices to apply $\exp_F$ to the equation $\text{Tr}_{K_{n+1}/K_n}(a_{n+1}) = \frac{q}{\pi}a_n$ noting that all of the terms $a_i$ satisfy $|a_i|<r$. Again from the discussion in 2.1 we have that $\exp_F$ takes addition to addition from the formal group law $F$. This implies that the sequence $(\alpha_n)$ will satisfy 

$$\text{Tr}^{\text{LT}}_{K_{n+1}/K_n}(\alpha_{n+1}) = [\frac{q}{\pi}](\alpha_n)$$

\noindent
hence $(\alpha_n)\in S$.
\\
\\
Now suppose $\alpha = (\alpha_i) = (\exp_F(a_i))$ can be interpolated so that there exists some series $f(x) \in \mathcal{O}_K[[x]]$ such that $f(u_i) = \alpha_i$ where $u_i$ is a compatible sequence of torsion points of $F$. We have $|\alpha_i| = |a_i|$ because all $a_i$ satisfy $|a_i| < r$. Then $|\alpha_i| = |a_i| \le |\frac{q}{\pi^2}|^i|a_0| = |\frac{q}{\pi^2}|^i|\alpha_0|$ so it is clear that $|f(u_i)| = |\alpha_i| \rightarrow 0 $ as $i \rightarrow \infty$ in the case that $\pi^3\mid q$. However by lemma 3.2.1 we must also have $\lim_{i\rightarrow \infty}|\alpha_i|$ exists and is nonzero if $(\alpha_i)$ is interpolated. This is a contradiction, so we get that such $(\alpha_i)$ cannot be interpolated.

\pagebreak

\subsection{Mapping Coleman series into the kernel of $\mathscr{L}$}

In this section we look at the $\mathbb{Z}_p$-module of power series interpolating norm compatible sequences of principal units in a tower of Lubin-Tate extensions, and we show that this module modulo certain series is isomorphic to a submodule of the kernel of Coleman's norm operator. Under the condition $q$ is not an integral power of $\pi$, so there does not exist $n \in \mathbb{Z}$ such that $q = \pi^n$, we are able to show there is an injection from the $\mathbb{Z}_p$-module of norm compatible sequences of principal units into the kernel of Coleman's trace operator.
\\
\\
Let $\mathscr{A}_{G_m}$ be the $\mathbb{Z}_p$-module of all series $f \in \mathcal{O}_K[[x]]$ satisfying there exists a norm compatible sequence of principal units $(x_n)_{n \ge 0}$ with $x_n \in K_n$ such that $f(u_n) = x_n$ for all $n \ge 0$. In this section we will just denote $\mathscr{A}_{G_m}$ by $\mathscr{A}$ when there is no risk of confusing it with the module with the same name from section 2.1. We define a map from $\mathscr{A}_{G_m}$ to the kernel of Coleman's norm operator, and we show that the kernel of this map is either trivial or generated by a single series as a $\mathbb{Z}_p$-module.
\\
\\
Note that in order for $\mathscr{A}$ to be considered as a $\mathbb{Z}_p$-module we must define the action of $\mathbb{Z}_p$ on power series $f \in \mathcal{O}_K[[x]]$ satisfying $f(0) \equiv 1 \mod \pi$. It suffices to check that the sequence of series $f(x)^{p^n}$ converges coefficientwise to $1$, and then apply lemma 5.2 to guarantee the coefficientwise convergence of $f(x)^\alpha$ where $\alpha \in \mathbb{Z}_p$. To check that the sequence $(f(x)^{p^n})$ converges to $1$ it suffices to apply lemma 5.2. Let $m$ be any element of $\pi\mathcal{O}_K$ so that $f(m)$ is a principal unit. The sequence $(f(m)^{p^n})$ converges to $1$, so we get that the sequence of series $(f(x)^{p^n})$ also converges coefficientwise to $1$ by the lemma. It follows that $f(x)^\alpha$ is a well-defined power series in $\mathcal{O}_K[[x]]$ whenever $\alpha \in \mathbb{Z}_p$ and $f(0) \equiv 1 \mod \pi$. From the above we also get $\mathbb{Z}_p$ acts on series $g \in \mathcal{O}_K[[x]]$ satisfying $\pi \mid g(0)$ by the following definition:

$$[\alpha](g(x)) = (1 + g(x))^\alpha - 1$$

\noindent
for every $\alpha \in \mathbb{Z}_p$.
\\
\\
As in previous sections fix a uniformizer $\pi$ of $\mathcal{O}_K$. Fix some choice of series $f_0(x) \in \mathcal{O}_K[[x]]$ such that $f_0(x) \equiv x^q \mod \pi$ and such that $f_0(x) \equiv \pi x \mod \deg 2$. $F$ will always denote the Lubin-Tate formal group law associated to $f_0(x)$. $x\oplus_F y$ will always denote the operation $F(x,y)$, the subscript denoting addition coming from the formal group law $F$. $x\oplus_{G_m} y$ will denote addition coming from the multiplicative formal group law. Any $\oplus$ appearing without a subscript will denote addition with respect to the multiplicative formal group law for this section. 
\\
\\
For this section $[a]_F$ always denotes the element of End$(F)$ corresponding to $a\in \mathcal{O}_K$, and $[a]_{G_m}$ will always denote the element of End$(G_m)$ corresponding to $a\in \mathbb{Z}_p$. For this section only we interpret $[a]$ to mean the endomorphism of $G_m$ corresponding to $a \in \mathbb{Z}_p$ if there is no subscript. $\textfrak{F}_n$ will always denote the $n$-th level torsion points of $F$, whereas if we need to write down the $n$-th level torsion points of $G_m$ we will write $(\textfrak{G}_m)_n$.
\\
\\
We must study the map $\phi_{G_m}(f) = [q]_{G_m}(f(x))\ominus_{G_m} f([\pi]_F(x))$ for $f$ where $1+f(x) \in \mathscr{A}_{G_m}$ and find its kernel. We prove the following:

\begin{theorem*}
The map $\mathscr{A}_{G_m} \rightarrow \text{ker}(\mathscr{L}) \subset \mathcal{O}_K[[x]]$ defined by sending $g \in \mathscr{A}_{G_m}$ to $\log([p^r]\phi_{G_m}(g - 1))$ is an injection of the $\mathbb{Z}_p$-module of all norm compatible sequences of principal units into the kernel of Coleman's trace operator when $q$ is not an integer power of $\pi$. The kernel of the same map is either trivial or generated by a single series as a $\mathbb{Z}_p$-module if $q = \pi^n$ for some $n$.

\end{theorem*}

\noindent
In the above theorem $r$ is some fixed integer such that 

$$[p^r](\pi\mathcal{O}_K[[x]]) \subseteq p\mathcal{O}_K[[x]]$$

\noindent
and $\log$ denotes $\log_{G_m}$ the logarithm of the multiplicative formal group law.
\\
\\
$(\textfrak{G}_m)_\infty$ denotes the set of all torsion points of the multiplicative formal group law over $\mathbb{Q}_p$. Let $\mathscr{E}(G_m)$ denote the set of all sequences $(a_i)_{i \ge 0}$ such that $a_i \in (\textfrak{G}_m)_\infty$ for each $i$, $[q](a_{i+1}) = a_i$ for each $i$, and $[q](a_0) = 0$. We show that if $|f(0)| < 1$ and $[q]_{G_m}(f(x))\ominus_{G_m} f([\pi]_F(x)) =0 $ then $f$ interpolates some sequence in $\mathscr{E}(G_m)$. In particular the set of all such $f$ in the kernel is either empty or generated by a single series as a $\mathbb{Z}_p$-module.
\\
\\
Suppose $f$ is in the kernel of $\phi_{G_m}$, so that

$$[q](f(x))\ominus f([\pi]_F(x)) = 0$$ 

\noindent
Then $[q-1]_{G_m}(f(0)) = 0$ which is only possible if $f(0) = 0$ since $[q-1]_{G_m}$ is an isomorphism on the maximal ideal in $\Omega$. This implies $[q]_{G_m}(f(u_0)) = 0$, which is only possible if $f(u_0) = a_0 \in (\textfrak{G}_m)_\infty$ and $[q]_{G_m}(a_0) = 0$. Since $f$ is in the kernel of the above map we get the recursive relationship $[q]_{G_m}(f(u_{i+1})) = f(u_i)$, so we know that the sequence $(f(u_i))$ is some sequence in $\mathscr{E}(G_m)$. 
\\
\\
For all $n$ we label the $n$-th level torsion points of $G_m$ by $(\textfrak{G}_m)_n$. $\mathbb{Z}_p^\times$ acts transitively on all sequences of the form $(a_i)_{i \ge 0}$ where $a_n \in (\textfrak{G}_m)_n$ and the $a_n$ satisfies $[p](a_{n+1}) = a_n$. Let $f$ be a series in the kernel of $[q]_{G_m}(f(x))\ominus_{G_m} f([\pi]_F(x))$, so that $f$ must interpolate some sequence in $\mathscr{E}(G_m)$. We will write the index of $f$ to denote the smallest $n$ such that $f(u_0) \in (\textfrak{G}_m)_n$. Since we know $[q]_{G_m}(f(u_0)) = 0$, we know there is an upper bound for the index of $f$ over all choices of $f$. Pick $f_0$ to be any such $f$ with maximal index. We show for any $f$ in the kernel there exists some $a \in \mathbb{Z}_p$ such that $[a]_{G_m}(f_0) = f$.
\\
\\
Let $i_f$ denote the index of $f$. Since $f_0$ was chosen to have maximal index we know there exists some power of $p$, say $p^t$, such that $[p^t]_{G_m}(f_0(u_0))$ lives in $(\textfrak{G}_m)_{i_f}$. Then because $[p]((\textfrak{G}_m)_{n+1}\setminus (\textfrak{G}_m)_n) = (\textfrak{G}_m)_n \setminus (\textfrak{G}_m)_{n-1}$ for all $n \ge 1$ it follows for arbitrary $n \ge 0$ we have $[p^t]_{G_m}(f_0(u_n))$ and $f(u_n)$ will have the same level as torsion points of $G_m$. If $i_f > 0$ we shift both sequences by the endomorphism $[p^{i_f}]$, so we consider the sequences $[p^{t+i_f}]_{G_m}(f_0(u_n))$ and $[p^{i_f}]_{G_m}f(u_n)$. For $n=0$ both entries live in $(\textfrak{G}_m)_0$.
\\
\\
Now consider the set of all sequences $(b_n)_{n \ge 0}$ such that $b_n \in (\textfrak{G}_m)_n$, and $[p](b_{n+1}) = b_n$ for all $n$. Let $(\textfrak{G}_m)_\infty$ denote the set of all torsion points of $G_m$. $\mathbb{Q}_p( (\textfrak{G}_m)_\infty )$ is an abelian extension of $\mathbb{Q}_p$ with Galois group isomorphic to $\mathbb{Z}_p^\times$ by local class field theory. This Galois group acts transitively on the set of all such sequences. Since the sequences $[p^{t+i_f}](f_0(u_n))$ and $[p^{i_f}]f(u_n)$ are both subsequences of such sequences, and $[p^{t+i_f}](f_0(u_n))$ and $[p^{i_f}]f(u_n)$ both have the same level as torsion points of $G_m$, we get that there must exist some $u \in \mathbb{Z}_p^\times$ such that $[up^{t+i_f}](f_0(u_n)) = [p^{i_f}]f(u_n)$ for all $n$. Then since the series $[up^{t+i_f}](f_0(x))$ and $[p^{i_f}](f(x))$ agree on all torsion points of $F$, they must be equal in $\mathcal{O}_K[[x]]$. This is only possible if $[up^t](f_0(x)) = f(x)$ since $[p](x)$ has a formal power series inverse under composition.
\\
\\
The above completes the proof that the $\mathbb{Z}_p$-module of series in $\mathcal{O}_K[[x]]$ which interpolate sequences in $\mathscr{E}(G_m)$ is either empty or it is generated by a single series.
\\
\\
For the remainder of this section we let $\mathscr{A}'$ denote the set of all power series $f(x) \in \mathcal{O}_K[[x]]$ satisfying $|f(0)| < 1$ and satisfying the identity

$$f(x) \oplus_{G_m} f(x\oplus_F z_1) \oplus_{G_m} \ldots \oplus_{G_m} f(x \oplus_F z_{q-1}) = f([\pi]_F(x))$$

\noindent
the $z_i$ in the above identity range over all elements of $\textfrak{F}_0$.
\\
\\
Note that $f \in \mathscr{A}'$ is equivalent to the series $g(x) = 1 + f(x)$ interpolating some norm compatible sequence of principal units in the tower of field extensions $K_n = K(\textfrak{F}_n)$ over $K$. Conversely if $g(x)$ is a series in $\mathcal{O}_K[[x]]$ interpolating a norm compatible sequence of principal units then $f(x) = g(x) - 1$ must live in $\mathscr{A}'$.
\\
\\
Let $\mathscr{D}_{G_m,F}$ denote the set of all series $h(x) \in \mathcal{O}_K[[x]]$ satisfying the identity

$$\sum^{G_m}_{z\in \textfrak{F}_0} h(x \oplus_F z) = 0$$

\noindent
In this identity $\sum\limits^{G_m}$ denotes addition using the multiplicative formal group law.
\\
\\
For what follows we will need the definition of the norm operator from \cite{Coleman1}. 

\begin{theorem*}
There exists a unique map $\mathscr{N} : \mathcal{O}_K((x)) \rightarrow \mathcal{O}_K((x))$ which satisfies

$$\mathscr{N}(f)_{\pi} = \prod_{z\in\textfrak{F}_0} f(x\oplus_F z)$$

\noindent
Moreover, $\mathscr{N}$ is continuous.
\end{theorem*}

\noindent
Here the notation $f_\pi$ denotes the composition $f_\pi = f([\pi](x))$ for arbitrary $f \in \mathcal{O}_K((x))$. Note that $\mathscr{D}_{G_m, F}$ is isomorphic to the kernel of Coleman's norm operator in $\mathcal{O}_K[[x]]$. If $h \in \mathscr{D}_{G_m, F}$ then the series $h_0(x) = 1+ h(x)$ satisfies 

$$\mathscr{N}(h_0)_{\pi} = \prod_{z\in\textfrak{F}_0} h_0(x \oplus_F z) = \prod_{z \in \textfrak{F}_0}(1 + h(x \oplus_F z)) = (1 + \sum^{G_m}_{z\in\textfrak{F}_0} h(x \oplus_F z) ) = 1$$

\noindent
Similarly if $h_0$ is a series in the kernel of the norm operator, then $h(x) = -1 +h_0(x)$ satisfies the identity 

$$ \prod_{z \in \textfrak{F}_0} (h(x \oplus_F z) + 1) = 1$$

\noindent
hence

$$\sum^{G_m}_{z\in\textfrak{F}_0} h(x\oplus_F z) = 0$$

\noindent
We define $\phi_{G_m} : \mathscr{A}' \rightarrow \mathscr{D}_{G_m, F}$ to be the map:

$$\phi_{G_m}(f) = [q]_{G_m}(f(x)) \ominus_{G_m} f([\pi]_F(x))$$

\noindent
We must check that the image of $\phi_{G_m}$ is actually contained in $\mathscr{D}_{G_m, F}$. We have $\phi_{G_m}(f)(x \oplus_F z) = [q]_{G_m}(f(x \oplus_F z)) \ominus_{G_m} f([\pi]_F(x))$. This implies the series

$$ \sum^{G_m}_{z\in \textfrak{F}_0}\phi_{G_m}(f)(x \oplus_F z)$$

\noindent
must equal

$$ [q]_{G_m}(\sum^{G_m}_{z\in\textfrak{F}_0} f(x\oplus_F z)) \ominus_{G_m} [q]_{G_m}(f([\pi]_F(x)))$$

\noindent
Since the above is just $[q]_{G_m}$ applied to the equation defining $\mathscr{A}'$ we must have that it equals zero, hence $\phi_{G_m}(\mathscr{A}') \subseteq \mathscr{D}_{G_m, F}$.
\\
\\
One can check that $\phi_{G_m}$ is a map of $\mathbb{Z}_p$-modules. We already saw that the kernel of $\phi_{G_m}$ on $\mathscr{A}'$ was the submodule of all series interpolating a sequence from $\mathscr{E}(G_m)$. We label the set of series in $\mathscr{A}'$ which interpolate any element in $\mathscr{E}(G_m)$ by $\mathscr{E}$, so that $\mathscr{E}$ is either empty or generated by a single series as a $\mathbb{Z}_p$-submodule of $\mathscr{A}'$. 
\\
\\
We currently have that $\mathscr{A}'/\mathscr{E}$ is isomorphic to some submodule of $\mathscr{D}_{G_m,F}$ using the map $\phi_{G_m}$. In cases where we can show $\mathscr{E}$ is empty, we will use the same map to obtain an injection of $\mathscr{A}'$ into the kernel of Coleman's trace operator.
\\
\\
For the remainder of this section we consider the $\mathbb{Z}_p$-module of norm compatible sequences of principal units contained in the tower of field extensions $K_n = K(\textfrak{F}_n)$. We denote the $\mathbb{Z}_p$-module of all such series interpolating any such norm compatible sequence as $\mathscr{A}_{G_m}$ or just $\mathscr{A}$ for short when there is no risk of confusion in the remainder of this section. Note that if $f \in \mathscr{A}$ then the series given by $-1 + f(x)$ interpolates some sequence of the form $\alpha_i \in K_i$ where 

$$\sum^{G_m}_{g \in G(K_{i+1}/K_i)} g(\alpha_{i+1}) = \alpha_i$$

\noindent
and each $\alpha_i$ lives in the maximal ideal of $\mathcal{O}_{K_i}$. Conversely if $f$ interpolates a sequence $(\alpha_i)$ of the above form (so that $f(u_i) = \alpha_i$), we get that the series $1 + f(x)$ interpolates some norm compatible sequence of principal units. We use $\mathscr{A}'$ to denote the $\mathbb{Z}_p$-module of all series $\mathscr{A}' = \{f \mid 1+f(x) \in \mathscr{A} \}$. This definition of $\mathscr{A}'$ is equivalent to the definition previously given in this section. For all $f$ living in $\mathscr{A}'$ we must have $|f(0)| < 1$ becaue $|f(u_i)| < 1$ for all $i$ and because $f \in \mathcal{O}_K[[x]]$.
\\
\\
We refer to the map $\phi_{G_m}$ as $\phi$ when there is no risk of confusion between maps with the same name from other sections. We prove that $\phi$ is an injection under the condition $\pi^n \neq q$ for all positive integer exponents $n$. In this case suppose

$$f(x) = \sum_{n=0}^\infty a_nx^n$$

\noindent
is such that $\phi(f) = 0$. We show that $a_n = 0$ for each $n$. First note that $\phi(f)(0) = [q](a_0)\ominus a_0$ can only be zero if $a_0 = 0$. This is because we have the condition $|a_0| < 1$. Now suppose we have shown $a_i = 0$ for each $i$ with $0 \le i \le N$ for some nonnegative integer $N$. We show this implies $a_{N+1}$ must also be zero.
\\
\\
Consider the expression $\phi(f)$ mod deg $N+2$. Since all terms of $f$ divisible by $x^{N+2}$ will vanish when we consider the expression mod deg $N+2$, we get that

$$\phi(f)(x) \equiv [q]_{G_m}( \sum_{n=0}^{N+1} a_nx^n ) \ominus_{G_m} ( \sum_{n=0}^{N+1} a_n[\pi]_F(x)^n ) \text{ mod deg } N+2$$

\noindent
Since we have that $a_i = 0$ for $0 \le i \le N$ by our inductive hypothesis the above can be rewritten as 

$$\phi(f)(x) \equiv [q]_{G_m}( a_{N+1}x^{N+1} ) \ominus_{G_m} (a_{N+1}[\pi]_F(x)^{N+1}) \text{ mod deg } N+2$$

\noindent
Now we note $[q]_{G_m}(a_{N+1}x^{N+1}) \equiv qa_{N+1}x^{N+1} \text{ mod deg } N+2$ and also that $a_{N+1}[\pi]_F(x)^{N+1} \equiv a_{N+1}\pi^{N+1}x^{N+1} \text{ mod deg } N+2$. Substituting these terms into the previous expression for $\phi(f)$ gives:

$$ \phi(f)(x) \equiv (qa_{N+1}x^{N+1}) \ominus_{G_m} (a_{N+1}\pi^{N+1}x^{N+1}) \text{ mod deg } N+2$$

\noindent
Now the nonlinear term in the expansion of the multiplicative formal group law above vanishes mod deg $N+2$, so we get that

$$\phi(f)(x) \equiv (qa_{N+1}x^{N+1}) + i_{G_m}(a_{N+1}\pi^{N+1}x^{N+1}) \text{ mod deg } N+2$$

\noindent
where $i_{G_m}$ is the unique series for the multiplicative formal group law satisfying 

$$ x + i_{G_m}(x) + xi_{G_m}(x) = 0 $$

\noindent
for all $x$. Since $i_{G_m}(x) \equiv -x \mod x^2$ we are able to write  

$$\phi(f)(x) \equiv (qa_{N+1}x^{N+1}) - (a_{N+1}\pi^{N+1}x^{N+1}) \text{ mod deg } N+2$$

\noindent
Therefore we have 

$$\phi(f)(x) \equiv a_{N+1}(q-\pi^{N+1})x^{N+1} \text{ mod deg } N+2$$

\noindent
We are working under the assumption $q-\pi^{N+1} \neq 0$, hence $\phi(f)(x)$ is not congruent to zero mod deg $N+2$ if $a_{N+1} \neq 0$. Since this contradicts that $\phi(f) = 0$, we must have that $a_{N+1} = 0$. This completes the proof that the kernel of $\phi$ is trivial when $\pi$ is any uniformizer satisfying $\pi^n \neq q$ for all $n$. If we are in the case where $\pi^n = q$ for some positive integer $n$ then kernel of $\phi$ is still either empty or generated by a single series as shown earlier in this section.
\\
\\
Next we show that for all choices of $f \in \mathscr{A}'$ we have $\pi \mid \phi(f)$ in $\mathcal{O}_K[[x]]$. First consider $[q]_{G_m}(f(x)) \mod \pi$. Because $[p]_{G_m}(x) \equiv x^p \mod p$ and $[p]_{G_m}(x) \equiv px \mod \deg 2$ we get that $[q]_{G_m}(x) \equiv x^q \mod \pi$. This implies $[q]_{G_m}(f(x)) \equiv f(x)^q \mod \pi$. We also have that $[\pi]_F(x) \equiv x^q \mod \pi$. This implies that $f([\pi]_F(x)) \equiv f(x^q) \mod \pi$. Putting these together gives

$$ \phi(f)(x) = [q]_{G_m}(f(x)) \ominus_{G_m} f([\pi]_F(x)) \equiv f(x)^q \ominus_{G_m} f(x^q) \mod \pi$$

\noindent
Then we have $f(x)^q \equiv f(x^q) \mod \pi$ because $q$ is the size of the residue field $\mathcal{O}_K/\pi\mathcal{O}_K$ and $f(x) \in \mathcal{O}_K[[x]]$. We conclude that $\phi(f) \equiv 0 \mod \pi$.
\\
\\
Now we note that if $\pi^n \mid f(x)$ in $\mathcal{O}_K[[x]]$ and $n \ge 1$ we must have $\pi^{n+1} \mid [p]_{G_m}(f(x))$. This is true because for $n \neq p$ we have $p$ divides the coefficient of $x^n$ appearing in $[p]_{G_m}(x)$. For the term $n = p$ we have $\pi^{np} \mid (f(x))^p$ if $\pi^n \mid f(x)$. This is enough to show that $\pi^{n+1} \mid [p]_{G_m}(f(x))$ whenever $\pi^n \mid f(x)$. From the above it follows that there exists a positive integer $r$ such that

$$[p^r]_{G_m}(\phi(f)(x)) \in p\mathcal{O}_K[[x]]$$

\noindent
for all choices of $f \in \mathscr{A}'$. For example if $r$ is any integer sufficiently large so that $p \mid \pi^r$ then the above argument shows that 

$$ \pi^r \mid [p^r]_{G_m}(\phi(f)(x))$$

\noindent
hence $[p^r]_{G_m}(\phi(f)(x)) \in p\mathcal{O}_K[[x]]$ for all $f \in \mathscr{A}'$ for such a choice of $r$. From this point on we take $r$ to be the unique smallest integer satisfying the above condition.
\\
\\
We will denote $\mathscr{D}_{G_m}$ by $\mathscr{D}$ when there is no risk of confusion with modules by the same name in other sections. Note that $\phi(f) \in \mathscr{D}$ implies $[p^r]_{G_m}(\phi(f)) \in \mathscr{D}$ as well. This is because

$$ \sum^{G_m}_{z\in\textfrak{F}_0} [p^r]_{G_m}(\phi(f)(x\oplus_F z)) = [p^r]_{G_m}( \sum^{G_m}_{z\in\textfrak{F}_0} \phi(f)(x\oplus_F z)) = [p^r]_{G_m}(0) = 0$$

\noindent 
From this point on we will always denote the logarithm of the multiplicative formal group law by $\log(x) = \log_{G_m}(x)$. We show that $\log : \mathscr{D}\cap p\mathcal{O}_K[[x]] \rightarrow M'$ is an ismorphism of $\mathbb{Z}_p$-modules. Here we use $M'$ to denote the $\mathbb{Z}_p$-module of all series in $g \in p\mathcal{O}_K[[x]]$ satisfying $\mathscr{L}_F(g) = 0$. This is the same as writing $M'$ is the $\mathbb{Z}_p$-module of all series $g \in p\mathcal{O}_K[[x]]$ satisfying

$$ \sum_{z\in\textfrak{F}_0} g(x\oplus_F z) = 0$$

\noindent
Now since $p \mid h(x)$ we have that the composition of series $\log(f(x)) \in p\mathcal{O}_K[[x]]$. This is because it is well known that $\log_{G_m} : p \mathcal{O}_K \rightarrow p\mathcal{O}_K$ and $\exp_{G_m} : p \mathcal{O}_K \rightarrow p\mathcal{O}_K$ are inverse isomorphisms of $\mathbb{Z}_p$-modules. See for example Proposition 7.17 and Proposition 2.4 in \cite{Kolyvagin}, agreeing with the $\mathbb{Z}_p$-action follows if we consider Theorem 2 in section 5.1 of \cite{BorevichShafarevich}. The same estimates of divisibility of $\log_F(a)$ and $\exp_F(b)$ depending on divisibility of $a,b$ in the above proof imply that $\log = \log_{G_m} : p\mathcal{O}_K[[x]] \rightarrow p\mathcal{O}_K[[x]]$ and $\exp = \exp_{G_m} : p\mathcal{O}_K[[x]] \rightarrow p\mathcal{O}_K[[x]]$ are defined (as coefficientwise limits, see the beginning of section 2.1). The remaining claims follow because they are free for substitutions $x\in p\mathcal{O}_K$, and coefficientwise limits agree with composition of functions on $p\mathcal{O}_K$, and $\log_F(\mathscr{D}\cap p \mathcal{O}_K[[x]]) = M'$. 
\\
\\
Choose a uniformizer $\pi$ satisfying $\pi^n \neq q$ for all positive integers $n$. Then the conclusion in this case is the following:

\begin{theorem}
If $q \neq \pi^n$ for all integers $n$ then the map taking $f \in \mathscr{A}'$ to $\log( [p^r](\phi(f)(x)))$ is an injection of the $\mathbb{Z}_p$-module of all series interpolating norm compatible sequences of principal units into $M'$, the $\mathbb{Z}_p$-module of all series in the kernel of Coleman's trace operator and divisible by $p$.
\end{theorem}

\noindent
If we again compose the above map with one of the maps defined by $t_n : M' \rightarrow M'$ defined by $t_n(g) = ([\pi](x))^ng(x)$ for $n \ge 2$ we obtain an injection of $\mathscr{A}'$ into the $\mathcal{O}_K$-module $\mathscr{C}'$, defined in section 2.1 to be the $\mathcal{O}_K$-module of all series $g \in \mathcal{O}_K[[x]]$ satisfying $\mathscr{L}_F(g) = 0$ and $g'(0) = 0$. This shows that the set of interpolated sequences in $S$ is at least as large as the set of all norm compatible sequences of principal units over the same tower of field extensions at least when $\pi$ satisfies the previously mentioned condition.

\pagebreak

\numberwithin{theorem}{section}

\section{Joint Work with Victor Kolyvagin}

\subsection{Explicit interpolation theorem}

Let $K$ be a finite extension of $\mathbb{Q}_p$, and let $\pi$ be a uniformizer of $K$. Let $\mathcal{O}_K$ denote the ring of integers in $K$. Let $q$ be the size of the residue field $\mathcal{O}_K/\pi\mathcal{O}_K$. Fix some series $f(x) \in \mathcal{O}_K[[x]]$ such that $f(x) \equiv \pi x \mod \deg 2$, and $f(x) \equiv x^q \mod \pi$. We will let $F_f = F_f(x,y)$ denote the Lubin-Tate formal group law associated to $f$. We will use $x \oplus_f y = F_f(x,y)$ to denote the formal group law operation. We let $\textfrak{F}_n$ denote the set of all $n$th level torsion points of of $F_f$, so that $\textfrak{F}_n$ is the set of all zeros of the series $f^{(n+1)}(x)$.
\\
\\
Let $R = R(\pi, q)$ be the set of all $g \in \mathcal{O}_K[x]$ such that $g(x) \equiv \pi x \mod \deg 2$, and such that $g$ is a monic polynomial of degree $q$ and $g(x) \equiv x^q \mod \pi$.

\begin{proposition}
    Let $r(x) \in \mathcal{O}_K[[x]]$ with $r(0) = 0$. Then

    $$\prod_{z \in \textfrak{G}_0}(r(x)\oplus_g z) = (-1)^{p-1}g(r(x))$$
\end{proposition}

\noindent
In the above there is a Lubin-Tate formal group law $F_g = F_g(x,y)$ associated to $g$, and $\oplus_g$ denotes the operation $x \oplus_g y = F_g(x,y)$. $\textfrak{G}_0$ denotes the zero level torsion points of $F_g$.
\\
\\
Proof of Proposition 4.1: it is enough to consider the case $r = x$, then the substitution $x \mapsto r$ proves the claim. The polynomial $g(Y) - g(X)$ considered as polynomial in $Y$ over the ring $\mathcal{O}_K[x]$ has roots $X \oplus_g z$ where $z$ spans over all elements of $z \in \textfrak{G}_0$. Hence we get

$$\prod_{z\in\textfrak{G}_0} -(X\oplus_g z) = -g(X)$$

\noindent
and the proposition follows because $|\textfrak{G}_0| = q = p^f$.
\\
\\
Let $r = i_{f,g}(x) \in \mathcal{O}_K[[x]]$ be the isomorphism of $F_f$ and $F_g$ such that $i_{f,g}(x) \equiv x \mod \deg 2$, so that $g(r) = r(f)$. See section 3.5 of chapter 6 of \cite{CasselsFrohlich} for the existence of such an isomorphism.

\begin{proposition}
    $$\prod_{v\in \textfrak{F}_0}r(x\oplus_f v) = (-1)^{p-1}r(f(x))$$

    \noindent
\end{proposition}

\noindent
Proof: 

$$\prod_{v \in \textfrak{F}_0} r(x\oplus_f v) = \prod_{v \in \textfrak{F}_0}(r(x) \oplus_g r(v)) = \prod_{z\in \textfrak{G}_0}(r(x)\oplus_g z)$$

\noindent
The above is equal to

$$\prod_{z\in \textfrak{G}_0}(r(x) \oplus_g z) = (-1)^{p-1}g(r(x)) = (-1)^{p-1}r(f(x))$$

\noindent
by Proposition 4.1. This concludes the proof of Proposition 4.2.
\\
\\
Now we also have

$$\prod_{v\in \textfrak{F}_0} (-1)^{p-1}r(x\oplus_f v) = (-1)^{(p-1)q} \prod_{v\in \textfrak{F}_0} r(x\oplus_f v) = (-1)^{p-1}r(f(x))$$

\noindent
It follows that if $v_n \in \textfrak{F}_n$ is any sequence such that $f(v_{n+1}) = v_n$ then the sequence $\{ (-1)^{p-1}r(v_n)\}$ for $n \ge 0$ is a norm compatible sequence such that $(-1)^{p-1}r(v_n) \in K_n = K(\textfrak{F}_n)$ for each $n$.
\\
\\
The series $(-1)^{p-1}i_{f,g}(x)$ where $g$ runs through the set $R(\pi, q)$ appear as "explicit" series with the property

$$\mathcal{N}_f(s) = s$$

\noindent
where $s \in \mathcal{O}_K[[x]]$ and $\mathcal{N}_f(s) \in \mathcal{O}_K[[x]]$ is defined by

$$\mathcal{N}_f(s)(f(x)) = \prod_{v\in \textfrak{F}_0} s(x\oplus_f v)$$

\noindent
These series provide a supply of "explicit" norm compatible systems in the tower of fields $\{K_n\}$.
\\
\\
The next step in the proof is to show that certain norm compatible systems $(a_n)$ are generated by the systems $((-1)^{p-1}i_{f,g}(v_n))$, where $g$ ranges over the elements of $R(\pi, q)$, with the consequence that there exists a series $r(x)$ such that $\mathcal{N}_f(r)=r$ and $a_n = r(v_n)$.
\\
\\
We first prove that for any norm compatible sequence $(a_n)$ the entry $a_0 \in K_0^{\times}$ is in the subgroup of $K_0^{\times}$ generated by the elements $(-1)^{p-1}i_{f,g}(v_0)$ where $g$ runs through $R$ and $v_0$ runs through $\textfrak{F}_0$.
\\
\\
Let $\xi$ be a uniformizer of the field $K_0$ such that

$$N_{K_0/K}(\xi) = \pi$$

\noindent
We know such uniformizers of $K_0$ exist because $v_0$ is a uniformizer of $K_0$, and we have the minimal polynomial of $v_0$ over $K$ is equal to

$$\frac{f(x)}{x} = x^{q-1} + \ldots + \pi$$

\noindent
It follows that $N_{K_0/K}(v_0) = \pi$ and $v_0$ is such a uniformizer when $q$ is odd. If $q$ is even we take our uniformizer to be $-v_0$ instead. Then let $g$ be $x$ times the characteristic polynomial of $(-1)^{p-1}\xi$ relative to the extensions $K_0/K$. Then $g$ is equal to

$$g(x) = x \prod_{\sigma \in G(K_0/K)}(x - (-1)^{p-1}\sigma(\xi)) = x(x^{q-1} + \ldots + (-1)^{p(q-1)}N_{K_0/K}(\xi))$$

\noindent
so that $(-1)^{p-1}\xi$ is a root of

$$g(x) = x^q + \ldots + \pi x$$

\noindent
and we have $g(x) \equiv x^q \mod \pi$. It follows that $g \in R(\pi, q)$ and $(-1)^{p-1}\xi$ being a nonzero root of $g$ must equal some $z_{0} \in \textfrak{G}_0\backslash \{0\}$. Hence $\xi = (-1)^{p-1}i_{f,g}(v_{0})$ where $v_{0} = i^{-1}_{f,g}(z_{0}) \in \textfrak{F}_0\backslash \{0\}$.
\\
\\
In the above we proved that the set of all elements $(-1)^{p-1}i_{f,g}(v_0)$ where $v_0$ runs over $\textfrak{F}_0$ and $g$ runs over $R(\pi, q)$ contains the set of all elements $\xi \in K_0$ such that $N_{K_0/K}(\xi) = \pi$, so that $\xi$ is a uniformizer of $K_0$.
\\
\\
Now suppose $(a_n)$ is a norm compatible system. Then 

$$a_0 \in \cap_n N_{K_n/K_0}(K_n^\times) = \{ b \in K_0^\times \mid N_{K_0/K}(b) \in \pi^{\mathbb{Z}} \}$$

\noindent
by the property (class field theory applied to the fields $K_n$) that

$$N_{K_n/K_0}(K_n^\times) = \{ c \in K_0^\times \mid N_{K_0/K}(c) \in \pi^{\mathbb{Z}}(1 + \pi^{n+1}\mathcal{O}_K)\}.$$

\noindent
Let $\xi$ be a uniformizer of $K_0$ with $N_{K_0/K}(\xi) = \pi$. Let $N_{K_0/K}(a_0) = \pi^j$. Then $a_0 = \xi^ju$ where $u$ is a unit of $K_0$ with $N_{K_0/K}(u) = 1$, so that $a_0 = \xi^{j-1} \xi_1$, where $N_{K_0/K}(\xi_1) = \pi$.
\\
\\
The above proves that $a_0$ is contained in the subgroup of $K_0^\times$ generated by $(-1)^{p-1}i_{f,g}(v_0)$ where $g$ runs through $R(\pi, q)$ and $v_0$ runs through $\textfrak{F}_0$. The next step is to show that a similar property holds for the entries $a_n$. 
\\
\\
The idea we will use to determine that certain $a_n$ are generated by the correct values is to split the extension $K_n/K_0$ into a tower

$$K_n = H_N \supset H_{N-1} \supset \ldots \supset H_1 \supset H_0 = K_0$$

\noindent
where each $H_j/H_{j-1}$ is a cyclic extension of degree $p$. This is possible because

$$G(K_n/K_0) \cong (1 + \pi \mathcal{O}_K)/(1+\pi^{n+1}\mathcal{O}_K)$$

\noindent
which is an abelian group of order $q^n$. We also have that 

$$G(K_{n}/K_{n-1}) \cong (1 + \pi^n\mathcal{O}_K) / (1 + \pi^{n+1}\mathcal{O}_K) \cong \mathbb{F}_q \cong (\mathbb{Z}/p\mathbb{Z})^f$$ 

\noindent
which allows us to construct such field extensions $H_j$.
\\
\\
Now let $\Delta$ be the cyclic subgroup of order $p-1$ in

$$G(K_0/K) \cong U_K/(1 + \pi \mathcal{O}_K) \cong \mathbb{F}_q^\times \cong \mathbb{Z}/(q-1)\mathbb{Z}$$

\noindent
If $A$ is a $\mathbb{Z}_p[\Delta]$-module and $\psi : \Delta \rightarrow \mu_{p-1} \subset \mathbb{Z}_p$ is a homomorphism, then 

$$A^\psi = \{ a \in A \mid \delta(a) = \psi(\delta)a \text{ for all } \delta \in \Delta\}$$

\noindent
We review the well known result that $A$ has a decomposition into eigenspaces defined by certain idempotents. Specifically we review a proof that $A = \bigoplus\limits_\psi A^\psi$ and $A^\psi = e_\psi A$ where $e_\psi \in \mathbb{Z}_p[\Delta]$ and

$$e_\psi = \frac{1}{p-1} \sum_{g \in \Delta} \psi(g^{-1})g $$

\noindent
Let $X$ denote the set of all of the homomorphisms $\psi : \Delta \rightarrow \mu_{p-1}$. To show the above decomposition $A = \bigoplus\limits_{\psi \in X} A^\psi$ it suffices to show

$$ \sum_{\psi \in X} e_\psi = 1$$

\noindent
and $e_{\psi_1}e_{\psi_2} = 0$ whenever $\psi_1 \neq \psi_2$. We also show $e_\psi^2 = e_\psi$. Then using that $A^\psi = e_\psi A$ gives the decomposition. 
\\
\\
We first show the inclusion $e_\psi A \subseteq A^\psi$. Take an arbitrary element $e_\psi a \in e_\psi A$ and apply any $\delta \in \Delta$. We have

$$\delta(e_\psi a) = \left(\frac{1}{p-1} \sum_{g \in \Delta}\psi(g^{-1})g\delta\right) a$$

\noindent
Then we can rewrite the sum

$$\left(\frac{1}{p-1} \sum_{g \in \Delta}\psi(g^{-1})g\delta\right) = \frac{1}{p-1} \sum_{g\in \Delta} \psi(\delta g^{-1})g = \frac{\psi(\delta)}{p-1}\sum_{g\in \Delta}\psi(g^{-1})g$$

\noindent
Substituting this into the previous equality for $\delta(e_\psi a)$ gives

$$\delta(e_\psi a) = \left(\frac{\psi(\delta)}{p-1}\sum_{g\in \Delta}\psi(g^{-1})g\right)a = \psi(\delta)e_\psi a$$

\noindent
so that $e_\psi a \in A^\psi$. 
\\
\\
We must now show the reverse inclusion. Take any $a \in A^\psi$. We show $a = e_\psi a$ so that $a \in e_\psi A$. Consider $e_\psi$ applied to $a$:

$$e_\psi a = \left(\frac{1}{p-1}\sum_{g\in\Delta}\psi(g^{-1})g\right)a$$

\noindent
Now for each $g \in \Delta$ we have $ga = \psi(g)a$ because $a \in A^\psi$. We get

$$\left(\frac{1}{p-1}\sum_{g\in\Delta}\psi(g^{-1})g\right)a = \frac{1}{p-1} \sum_{g\in\Delta} \psi(g^{-1})\psi(g)a = \frac{1}{p-1}\sum_{g\in\Delta}\psi(1)a = a$$

\noindent
This concludes the proof that $a = e_\psi a$ whenever $a \in A^\psi$, and we also conclude $A^\psi \subseteq e_\psi A$, so that we have $A^\psi = e_\psi A$.
\\
\\
Next we show $\sum\limits_{\psi \in X} e_\psi = 1$. First we have

$$\sum_{\psi \in X} e_\psi = \sum_{\psi \in X} \frac{1}{p-1}\sum_{g\in \Delta} \psi(g^{-1})g = \frac{1}{p-1} \sum_{\psi \in X} \sum_{g\in\Delta} \psi(g^{-1})g$$

\noindent
The above equals

$$\frac{1}{p-1}\sum_{g\in\Delta}\sum_{\psi\in X}\psi(g^{-1})g$$

\noindent
Now we show if $g$ is not the identity then  

$$\sum_{\psi \in X}\psi(g^{-1})g = 0$$

\noindent
Suppose $g$ is not the identity. Then there exists some $\psi_0$ such that $\psi_0(g^{-1}) \neq 1$. We have

$$\psi_0(g^{-1}) \sum_{\psi \in X} \psi(g^{-1}) = \sum_{\psi \in X} \psi_0(g^{-1})\psi(g^{-1}) = \sum_{\psi \in X} \psi(g^{-1})$$

\noindent
which is only possible if $\sum\limits_{\psi \in X} \psi(g^{-1}) = 0$. If $e \in \Delta$ is the identity, we conclude that

$$\sum_{\psi \in X}e_\psi = \frac{1}{p-1} \sum_{\psi \in X}\psi(1)e = 1$$

\noindent
We now move to the proof that $e_{\psi_1}e_{\psi_2} = 0$ whenever $\psi_1 \neq \psi_2$. We have

$$e_{\psi_1}e_{\psi_2} = \left(\frac{1}{p-1}\sum_{g \in \Delta} \psi_1(g^{-1})g\right) \left( \frac{1}{p-1}\sum_{h \in \Delta} \psi_2(h^{-1})h\right)$$

\noindent
The above is equal to

$$e_{\psi_1}e_{\psi_2} = \left(\frac{1}{p-1}\right)^2 \left( \sum_{g\in \Delta} \sum_{h \in \Delta} \psi_1(g^{-1})\psi_2(h^{-1})gh\right)$$

\noindent
which can be rewritten as 

$$\left( \frac{1}{p-1} \right)^2 \sum_{j \in \Delta} \sum_{h\in\Delta} \psi_1(hj^{-1})\psi_2(h^{-1})j = \left(\frac{1}{p-1}\right)^2 \sum_{j \in \Delta} \psi_1(j^{-1}) \sum_{h\in\Delta}\psi_1^{-1}\psi_2(h^{-1})j$$

\noindent
Now $\psi_1^{-1}\psi_2$ is some nontrivial element of $X$, so it suffices to prove that for every $\psi \in X$ with $\psi \neq 1$ we have 

$$\sum_{g\in\Delta}\psi(g^{-1}) = 0$$

\noindent
To prove the above note that $\psi \neq 1$ so there exists some $g_0 \in \Delta$ such that $\psi(g_0^{-1}) \neq 1$. We then have

$$\psi(g_0^{-1}) \sum_{g \in \Delta}\psi(g^{-1}) = \sum_{g\in\Delta} \psi(g_0^{-1}) \psi(g^{-1}) = \sum_{g\in\Delta} \psi(g^{-1})$$

\noindent
which is only possible if $\sum\limits_{g\in\Delta} \psi(g^{-1}) = 0$. We conclude that $e_{\psi_1}e_{\psi_2} = 0$ whenever $\psi_1 \neq \psi_2$.
\\
\\
We move to the proof that $e_\psi^2 = e_\psi$. Note that

$$e_\psi^2 = \left( \frac{1}{p-1} \right)^2 \left( \sum_{g\in\Delta} \psi(g^{-1})g \right) \left( \sum_{h\in\Delta} \psi(h^{-1})h \right)$$

\noindent
The above equals

$$\left( \frac{1}{p-1} \right)^2 \left( \sum_{g\in\Delta} \psi(g^{-1})g \right) \left( \sum_{h\in\Delta} \psi(h^{-1})h \right) = \left( \frac{1}{p-1} \right)^2 \sum_{g\in\Delta} \sum_{h \in \Delta} \psi(g^{-1})\psi(h^{-1})gh$$

\noindent
It then suffices to note

$$\sum_{g\in\Delta}\sum_{h\in\Delta} \psi(g^{-1}h^{-1})gh = (p-1)e_\psi$$

\noindent
which is true because for each $j \in \Delta$ there are exactly $p-1$ ways to write $j$ as a product $j = gh$ with $g\in\Delta$ and $h\in\Delta$. This concludes the proof of the decomposition for $\mathbb{Z}_p[\Delta]$-modules which we will use.
\\
\\
Let $U_{t,1}$ denote the principal units of $H_t$. Note that $U_{t,1}$ is a $\mathbb{Z}_p[\Delta]$-module. For each $t$ we fix a generator of $G(H_t/H_{t-1})$ which we denote by $\gamma_t$. We will need the following:

\begin{lemma}
    If $a\in U_{t,1}^\psi$ ($\psi \neq 1$) with $N_{H_t/H_{t-1}}(a)=1$ then $a = b^{\gamma_t-1}$ for some $b \in H_t^\times$ by Hilbert's 90th theorem. Furthermore this is true for some $b \in U_{t,1}^\psi$ if $\psi\neq \psi_0 = 1$. 
\end{lemma}

\noindent
Proof: $H_t^\times/(H_t^\times)^{p^m}$ is a $\mathbb{Z}/p^m\mathbb{Z}[\Delta]$-module with corresponding decomposition into $\psi$-eigenspaces. If $N_{H_t/H_{t-1}}(a) = 1$ we get $a = b^{\gamma_t - 1}$ for some $b \in H_t^\times$ by the Hilbert 90th theorem.
\\
\\
Now $a = a^{e_\psi} = (b^{e_{\psi,m}})^{\gamma_t-1}x^{p^m}$ where $e_{\psi,m}\in \mathbb{Z}[\Delta]$, and $e_{\psi,m} \equiv e_\psi$ (mod $p^m$). Let $\nu_{H_t}$ be the additive valuation on $H_t^\times$ normalized so that $\nu_{H_t}(\xi_t) = 1$ where $\xi_t$ is a uniformizer of $H_t$. Then 

$$\nu_{H_t}(b^{e_{\psi,m}}) = \left(\sum\limits_{\delta \in \Delta}\psi(\delta)\right)\nu_{H_t}(b) \equiv 0 \mod p^m$$ 

\noindent
for $\psi \neq \psi_0$ so we can pick $b^{e_{\psi,m}}$ to be a principal unit in $U_{t,1}$.
\\
\\
We get that $a = u_m^{\gamma_t-1}x_m^{p^m}$ where $u_m\in U_{t,1}$, and this implies $x_m\in U_{t,1}$ as well. Applying $e_\psi$ once more to the equation $a = u_m^{\gamma_t-1}x_m^{p^m}$ we can assume $u_m, x_m\in U_{t,1}^\psi$. 
\\
\\
Now $N_{H_t/H_{t-1}}(x^{p^m}_m) = 1$ implies $N_{H_t/H_{t-1}}(x_m) \in H_{t-1}^*\bigcap \mu_{p^m} \subseteq \mu_{p^N} = H_{t-1}^*\bigcap \mu_{p^\infty}$. $N$ does not depend on $m$ so that $N_{H_t/H_{t-1}}(x^{p^N}_m) = 1$. Then we can take $m = N+1$, and we obtain $a = u_{N+1}^{\gamma_t-1}(x_{N+1}^{p^N})^p$, where $x_{N+1}^{p^N}\in \ker N_{H_n/H_{n-1}}$. Then $1 = (x_{N+1}^{p^N})^{\sum\limits_{i=0}^{p-1}\gamma_t^i}$ so that 

$$x_{N+1}^{p^{N+1}} = (x_{N+1}^{p^N})^{p-\sum\limits_{i=0}^{p-1}\gamma_t^i} = \prod\limits_{i=0}^{p-1}(x_{N+1}^{p^N})^{(1-\gamma_t^i)} = z^{(\gamma_t-1)}$$

\noindent
where $z\in U_{t,1}^\psi$. This proves the lemma.

\begin{proposition}
 For the statement of this proposition we let $f_K$ denote the residue degree of the field $K$, so that $q = p^{f_K}$. Let $\psi$ be one of the nontrivial characters $\psi : \Delta \rightarrow \mu_{p-1}$. Suppose $u_0^j = (i_{f,g_{j_1}}(v_0)/i_{f,g_{j_0}}(v_0))^{e_\psi}$ spanning $j=1,\ldots,m$ generate $U_{0,1}^\psi$ as a $\mathbb{Z}_p$-module. Here we take the series $g_{j_0}$ and $g_{j_1}$ to be suitably chosen polynomials in $R(\pi,q)$. For $t = nf_K$ let $A_t$ be the $\mathbb{Z}_p[G(H_t/K_0)]$-module generated by $u_t^j = (i_{f,g_{j_1}}(v_{n})/i_{f,g_{j_0}}(v_{n}))^{e_\psi}$. For any $t$ not divisible by $f_K$ we let $u_t^j = N_{K_n/H_t}(u_{nf_K}^j)$ for any $n$ such that $nf_K > t$, and we define $A_t$ to be the $\mathbb{Z}_p[G(H_t/K_0)]$-module generated by $u_t^j$. Then we have $A_t = U_{t,1}^\psi$.
\end{proposition}

\noindent
Proof: Induction on $t$. It is true for $t=0$ by the condition. Suppose it is true for $U_{t-1,1}^\psi$. Let $u=u_0\in U_{t,1}^\psi$. Let $b_{t-1} \in A_{t-1}$ be the element such that $N_{H_t/H_{t-1}}(u) = b_{t-1}$. Then because $N_{H_t/H_{t-1}}(u_t^j)=u_{t-1}^j$ there exists $u_1\in A_t$ such that $N_{H_t/H_{t-1}}(u_1)=b_{t-1}$. Hence $N_{H_t/H_{t-1}}(u/u_1)=1$ implying $u=u_1x_1^{(\gamma_t-1)}$ with $x_1 \in U_{t,1}^\psi$ by lemma 4.3.
\\
\\
We apply the same reasoning to $x_1$ to get $u=u_1u_2^{(\gamma_t-1)}x_2^{(\gamma_t-1)^2}$, and by induction there exists $u_1, u_2, \ldots u_N \in A_t$ and $x_1, \ldots, x_N \in U_{t,1}^\psi$ for every $N$ such that $u = u_1u_2^{(\gamma_t-1)}\ldots x_N^{(\gamma_t-1)^N}$. Let $u_k = \prod\limits_{j=1}^m(u_t^j)^{\alpha_k^j}$ where $\alpha_k^j\in\mathbb{Z}_p[G(H_t/K_0)]$. Then 

$$u = \left(\prod\limits_{j=1}^m(u_t^j)^{\sum\limits_{k=1}^N\alpha_k^j(\gamma_t-1)^k}\right)x_N^{(\gamma_t-1)^N}$$

\noindent
We get that $u = \lim_{N\rightarrow\infty} u = \prod\limits_{j=1}^m(u_t^j)^{\sum\limits_{k=1}^\infty\alpha_k^j(\gamma_t-1)^k}$ where $(\gamma_t-1)^k\rightarrow 0$ as $k\rightarrow \infty$ in $\mathbb{Z}_p[G(H_t/H_{t-1})]$, so the proposition is proved.
\\
\\
Let us prove that there exist $u_0^j$, $j=1,\ldots,m$ which generate $U_{0,1}^\psi$ as a $\mathbb{Z}_p$-module. It suffices to show $U_{0,1}$ is finitely generated as a $\mathbb{Z}_p$-module. We have $\log : U_{0,1} \rightarrow K_0$ which induces an isomorphism $U_{0,\kappa} = 1+\pi_0^\kappa\mathcal{O}_{K_0}$ to $\pi_0^\kappa\mathcal{O}_{K_0}$ where $\kappa = \frac{q-1}{p-1}+1$. Now $\mathcal{O}_{K_0}$ is a free $\mathbb{Z}_p$-module of rank $[K_0/\mathbb{Q}_p]$, so $\pi_0^\kappa\mathcal{O}_{K_0}$ is also a free $\mathbb{Z}_p$-module of rank $[K_0/\mathbb{Q}_p]$. In addition $U_{0,1}/U_{0,\kappa}$ is finite hence $U_{0,1}$ and $U_{0,1}^\psi$ are finitely generated.
\\
\\
Now take any norm compatible sequence $(a_n)_{n \ge 0}$ with $a_n \in K_n$ where each $a_n$ is a principal unit. For each $\psi \in X$ we have taking the norm $N_{K_n/K_{n-1}} : K_n \rightarrow K_{n-1}$ commutes with applying $e_\psi$:

$$N_{K_n/K_{n-1}}(e_\psi a_n) = e_\psi N_{K_n/K_{n-1}}(a_n) = e_\psi a_{n-1}$$

\noindent
For this reason we get that each $\psi \in X$ gives us a norm compatible sequence $(e_\psi a_n)_{n \ge 0}$ where each $e_\psi a_n \in U_{nf,1}^\psi$. For $\psi \neq 1$ we can check that

$$N_{K_0/K}(e_\psi u) = 1$$

\noindent
for any principal unit $u \in U_{0,1}$ and this implies that the elements of the form $(-1)^pi_{f,g}(v_0)$ generate $U_{0,1}^\psi$. Putting this together with Proposition 4.4 implies we can find for each $N$ a series $f_N \in \mathcal{O}_K[[x]]$ such that $f_N(v_n) = e_\psi a_n$ for each $n$ with $n \le N$.
\\
\\
If $\psi \neq 1$ and $u \in U_{0,1}$ we have

$$N_{K_0/K}(e_\psi u) = e_\psi N_{K_0/K}(u) = \frac{1}{p-1} \left(\sum_{g\in \Delta} \psi(g^{-1})g\right)N_{K_0/K}(u)$$

\noindent
Note that because $N_{K_0/K}(u) \in K$ each $g \in \Delta$ acts trivially on $N_{K_0/K}(u)$. Then the above is equal to

$$N_{K_0/K}(e_\psi u) = N_{K_0/K}(u)^{\frac{1}{p-1}\sum\limits_{g \in \Delta} \psi(\frac{1}{g})}$$

\noindent
The above exponent, $\frac{1}{p-1}\sum\limits_{g\in\Delta}\psi(\frac{1}{g})$, is zero whenever $\psi$ is nontrivial. This concludes the proof that $N_{K_0/K}(e_\psi u) = 1$ whenever $u \in U_{0,1}$ and $\psi \neq 1$.
\\
\\
We now move from series which interpolate a sequence at the finite level to interpolating the entire sequence with a single series. If for all $N$ there exist norm compatible series $r_N \in \mathcal{O}_K[[x]]$ with $r_N(v_n) = a_n$ for all $n \le N$ then the sequence $\{r_N\}$ converges coefficient wise to some $r\in\mathcal{O}_K[[x]]$, and for this $r$ one can prove $a_n = r(v_{n})$ for all $n$.
\\
\\
We need the following lemma:

\begin{lemma}
    Let $r_n(x)$ be a sequence of series $r_n \in \mathcal{O}_K[[x]]$ and let $I$ denote the maximal ideal of the ring of integers of the algebraic closure of $K$. Then $(r_n(x))$ converges coefficientwise if and only if there exists a sequence $(m_i)$ with each $m_i \in I$ such that $\lim_{i\rightarrow\infty}|m_i| = 1$ and for each $m_i$ the sequence $(r_n(m_i))_{n\ge 0}$ converges.
\end{lemma}

\noindent
Proof: one direction is clear from lemma 5.2. For the other direction we assume we have a sequence of $m_i \in I$ such that $\lim_{i \rightarrow \infty}|m_i| =1$ and for each $i$ the sequence $(r_n(m_i))_{n\ge 0}$ converges.
\\
\\
Suppose $(r_n(x))$ does not converge coefficientwise. Label the coefficients of each $r_n$ so that

$$r_n(x) = \sum_{j=0}^\infty c_{n,j}x^j$$

\noindent
There is at least one index $j$ such that the sequence $(c_{n,j})_{n\ge 0}$ does not converge. For any such index let $M_j$ denote 

$$M_j = \limsup_{n\rightarrow \infty} |c_{n+1,j} - c_{n,j}|$$

\noindent
Then for any $j$ for which $(c_{n,j})_{n\ge 0}$ does not converge we must have $M_j \neq 0$. Now let $M$ be the maximal value of $M_j$ taken over all $j$. We then denote by $j_0$ the smallest index $j$ such that $M_{j_0} = M$. Now pick an element $m_i \in I$ such that $|m_i|^{j_0} > |\pi|$. We show that for such $m_i$ the sequence $(r_n(m_i))_{n\ge 0}$ cannot converge.
\\
\\
Consider the sequence of differences

$$r_{n+1}(m_i) - r_n(m_i) = \sum_{j=0}^\infty (c_{n+1,j} - c_{n,j})m_i^j$$

\noindent
For $j < j_0$ we have that $|c_{n+1,j}-c_{n,j}| < M$ for large enough $n$, which implies $|c_{n+1,j}-c_{n,j}| \le |\pi|M$ since all of the coefficients live in $\mathcal{O}_K$. Then for such $n$ we have 

$$|c_{n+1,j_0}-c_{n,j_0}| = M$$

\noindent
for infinitely many values of $n$. For these values of $n$ we get

$$|(c_{n+1,j_0}-c_{n,j_0})m_i^{j_0}| = M|m_i|^{j_0} > M|\pi| \ge |(c_{n+1,j}-c_{n,j})m_i^j|$$

\noindent
whenever $j < j_0$. Now pick $j_1$ to be an exponent large enough so that $|m_i|^{j_1} < M|m_i|^{j_0}$. Then for indices $j \ge j_1$ we get 

$$|(c_{n+1,j}-c_{n,j})m_i^j| < M|m_i|^{j_0} = |(c_{n+1,j_0}-c_{n,j_0})m_i^{j_0}|$$

\noindent
for the same infinitely many values of $n$ from above. We are then left with comparing $|(c_{n+1,j_0}-c_{n,j_0})m_i^{j_0}|$ to the terms $|(c_{n+1,j}-c_{n,j})m_i^j|$ with $j_0 < j <j_1$. Because $M$ is maximal, we know that for each $j$ with $j_0 < j < j_1$ there are only finitely many $n$ such that $|c_{n+1,j}-c_{n,j}| > M$. Therefore if we pick $n$ large enough we get that $|c_{n+1,j}-c_{n,j}| \le M$ for all indices $j$ in the above range. It follows that for infinitely many $n$ we have

$$|(c_{n+1,j_0}-c_{n,j_0})m_i^{j_0}| = M|m_i|^{j_0} > M|m_i|^j \ge |(c_{n+1,j}-c_{n,j})m_i^j|$$

\noindent
We get that for such values of $n$, the term $(c_{n+1,j_0}-c_{n,j_0})m_i^{j_0}$ is strictly the largest term in the expansion of $r_{n+1}(m_i) - r_n(m_i)$. Then

$$|r_{n+1}(m_i) - r_n(m_i)| = M|m_i|^{j_0}$$

\noindent
infinitely often which contradicts that $(r_n(m_i))$ converges. We get this contradiction from assuming $r_n(x)$ does not converge coefficientwise, so we conclude that $r_n(x)$ must converge coefficientwise if there is such a sequence $m_i$.
\\
\\
We conclude from lemma 4.5 and lemma 5.2 that if $(a_n)$ is a norm compatible sequence and $(f_N(x))$ is a sequence of series in $\mathcal{O}_K[[x]]$ such that $f_N$ satisfies $f_N(v_n) = a_n$ for all $n \le N$ then the coefficientwise limit of $(f_N(x))$ exists. We have that the coefficientwise limit $f(x)$ satisfies $f(v_n) = a_n$ for all $n$.

\pagebreak

\numberwithin{theorem}{section}

\section{Basic Lemmas}

\begin{lemma}
    Let $f$ and $g$ be two series in $\mathcal{O}_K[[x]]$ such that $f(x_i) = g(x_i)$ for infinitely many $x_i$ with $|x_i| < 1$. Then $f(x) = g(x)$
\end{lemma}

\noindent
Proof: let $d(x) = f(x) - g(x)$. If $d$ is nonzero there exists a power of $\pi$ such that $d(x) = \pi^n\delta(x)$ where $\delta(x) \in \mathcal{O}_K[[x]]$ and not all of the coefficients of $\delta(x)$ are divisible by $\pi$. By Weierstrass preparation theorem there exists a distinguished polynomial $F(x)$ and a unit $u(x) \in \mathcal{O}_K[[x]]$ such that 

$$\delta(x) = u(x)F(x)$$

\noindent
$u(0)$ is a unit, so that $u(x_i) \neq 0$ for each $i$. It follows that $F(x_i) = 0$ for each $x_i$ which is impossible. The assumption that $d$ is nonzero must be false, and it follows that $f(x) = g(x)$ in $\mathcal{O}_K[[x]]$.

\begin{lemma}
    Let $(f_n(x))$ be a sequence of power series in $\mathcal{O}_K[[x]]$. Let $I$ denote the maximal ideal of the ring of integers of the algebraic closure of $K$. Then the following are equivalent:

\begin{enumerate}
    \item $(f_n(x))$ converges coefficientwise
    \item for every $m \in I$ the sequence $(f_n(m))$ converges
    \item there exists a sequence $M$ of nonzero elements of $I$ converging to zero such that for each $m \in M$ the sequence $(f_n(m))$ converges
\end{enumerate}

\noindent
Also, in cases 2 and 3 there exists a unique $f(x) \in \mathcal{O}_K[[x]]$ such that $\lim f_n(m) = f(m)$ for $m \in I$ and $m \in M$ respectively, and $f(x)$ is the coefficientwise $\lim f_n(x)$.
\end{lemma}

\noindent
Proof: first we show $1$ implies $2$. It suffices to show for any $\epsilon > 0$ there exists $N$ such that whenever $n > N$ we have $|f_{n+1}(m) - f_n(m)| < \epsilon$. For each $n$ label the coefficients of $f_n(x)$ by

$$f_n(x) = \sum_{k=0}^\infty a_{n,k}x^k$$

\noindent
Since $m \in I$ there exists some $k_0$ such that $|m|^{k_0} < \epsilon$. Since $(f_n(x))$ converges coefficientwise we can find some constant $N$ such that $|a_{n+1,k} - a_{n,k}| < \epsilon$ for all $k$ with $k < k_0$ and for all $n > N$. Then for all $n$ greater than this choice of $N$ we have

$$f_{n+1}(m) - f_n(m) = \sum_{k=0}^\infty (a_{n+1,k}-a_{n,k})m^k$$

\noindent
It follows that if $k < k_0$ in the above we have

$$|(a_{n+1,k}-a_{n,k})m^k| \le |a_{n+1,k}-a_{n,k}| < \epsilon$$

\noindent
and if $k \ge k_0$ we have

$$|(a_{n+1,k}-a_{n,k})m^k| \le |m^k| \le |m|^{k_0} < \epsilon$$

\noindent
putting these together implies $|f_{n+1}(m) - f_n(m)| < \epsilon$ whenever $n > N$. This concludes the proof that $1$ implies $2$.
\\
\\
$2$ clearly implies $3$, so it suffices to show $3$ implies $1$ in order to show all three statements are equivalent. Suppose we have a sequence $M$ of nonzero elements of $I$ converging to zero such that $(f_n(m))$ converges for each $m \in M$. We show this implies the coefficientwise convergence of $(f_n(x))$.
\\
\\
Suppose that the sequence of power series $(f_n(x))$ does not converge coefficientwise. If $a_{n,k}$ is the $k$th coefficient of $f_n(x)$ so that

$$f_n(x) = \sum_{k=0}^\infty a_{n,k}x^k$$

\noindent
then the above implies there exists some $k$ such that the sequence $(a_{n,k})_{n\ge 0}$ is not Cauchy. Now take $k_0$ to be the minimal $k$ such that the sequence $(a_{n,k_0})_{n\ge 0}$ is not Cauchy. Then there exists some $\epsilon > 0$ such that

$$|a_{n+1,k_0} - a_{n,k_0}| \ge \epsilon$$

\noindent
for infinitely many values of $n$. Take some $m \in M$ satisfying $|m| < \epsilon$. We will show this gives a contradiction by showing $(f_n(m))$ is not Cauchy under the assumption $(f_n(x))$ does not converge coefficientwise.
\\
\\
We show that

$$|f_{n+1}(m) - f_n(m)| \ge |m|^{k_0}\epsilon$$

\noindent
for infinitely many values of $n$. First note that each sequence $(a_{n,k})_{n\ge 0}$ with $k < k_0$ is Cauchy because $k_0$ was assumed to be minimal. It follows there exists some $N$ such that 

$$|a_{n+1,k} - a_{n,k}| < |m|^{k_0}\epsilon$$

\noindent
whenever $n > N$ and $k < k_0$. Now we consider the difference

$$f_{n+1}(m) - f_n(m) = \sum_{k=0}^\infty (a_{n+1,k}-a_{n,k})m^k$$

\noindent
for any value of $n$ such that $n > N$ and such that

$$|a_{n+1,k_0} - a_{n,k_0}| \ge \epsilon$$

\noindent
Note that if $k > k_0$ we have

$$|(a_{n+1,k}-a_{n,k})m^k| \le |m|^k \le |m|^{k_0}|m| < |m|^{k_0}\epsilon$$

\noindent
Also for $k < k_0$ we have

$$|(a_{n+1,k}-a_{n,k})m^k| \le |(a_{n+1,k}-a_{n,k})| 
 < |m|^{k_0}\epsilon$$

\noindent
because $n > N$. It follows that

$$|f_{n+1}(m) - f_n(m)| = |(a_{n+1,k_0}-a_{n,k_0})m^{k_0}| \ge |m|^{k_0}\epsilon$$

\noindent
because all other terms appearing have strictly smaller absolute value. Since this inequality holds for infinitely many values of $n$ we conclude that $(f_n(m))$ cannot be Cauchy if $(f_n(x))$ does not converge coefficientwise. This completes the proof that $3$ implies $1$.
\\
\\
Now if any of the three equivalent condition holds then 1 holds and we let $f(x)$ be the unique series such that $(f_n(x))$ converges coefficientwise to $f(x)$. We define the function $h(x)$ on $I$ by setting 

$$h(m) = \lim_{n\rightarrow \infty}f_n(m)$$

\noindent
We will be able to show $h(m) = f(m)$ for each $m \in I$. We let 

$$f_n(x) = \sum_{k=0}^\infty a_{n,k}x^k$$

\noindent
It is possible to check that for each $m \in I$, the sequence $(f_n(m))$ converges to $f(m)$. Given $\epsilon > 0$ pick an integer $k_0$ such that $|m|^{k_0} < \epsilon$. Then if 

$$f(x) = \sum_{k=0}^\infty a_kx^k$$

\noindent
we can find some $N$ such that whenever $n > N$ we have 

$$|a_k - a_{n,k}| < \epsilon$$

\noindent
for all $k < k_0$. It follows for such $n$ that

$$f(m) - f_n(m) = \sum_{k=0}^\infty (a_k - a_{n,k})m^k$$

\noindent
In the above sum we have

$$|(a_k - a_{n,k})m^k| \le |m|^k < \epsilon$$

\noindent
if $k \ge k_0$, and we have

$$|(a_k - a_{n,k})m^k| \le |a_k - a_{n,k}| < \epsilon$$

\noindent
if $k < k_0$, so that we must have $|f(m) - f_n(m)| < \epsilon$. This completes the proof that

$$\lim_{n \rightarrow \infty} f_n(m) = h(m) = f(m)$$

\noindent
for all $m \in I$. This also shows that $f(x)$, the coefficientwise limit of $(f_n(x))$, is uniquely determined by the property $\lim f_n(m) = f(m)$ for all $m \in I$. To check this just note that if $g(x)$ is any other series satisfying $g(m) = \lim f_n(m)$ for all $m \in I$, then we must have $g$ and $f$ agree on infinitely many points in $I$, so they must be equal by lemma 5.1. This completes the proof of lemma 5.2.
\\
\\

\begin{lemma}
    If $f \in \mathcal{O}_K[[x]]$, then
    $$\mathscr{L}_F^n(f) \equiv 0 \mod \pi^n\mathcal{O}_K[[x]]$$

\noindent
This is lemma 6 from \cite{Coleman1}.
\end{lemma}

\noindent
Proof: the lemma follows directly from the case $n=1$ and the $\mathcal{O}_K$-linearity of $\mathscr{L}_F$. Let $p_0$ be the prime ideal of $\mathcal{O}_{K_0}$. Since $z \in p_0$ for any $z\in \textfrak{F}_0$ we have $f(x\oplus z) \equiv f(x) \mod p_0$. It follows that

$$\mathscr{L}_F(f)([\pi])(x) = \sum_{z\in \textfrak{F}_0} f(x\oplus z) \equiv q f(x) \equiv 0 \mod p_0$$

\noindent
Because both sides of the above congruence live in $\mathcal{O}_K[[x]]$ we get that

$$\mathscr{L}_F(f)([\pi](x)) \equiv 0 \mod \pi$$

\noindent
implying

$$\mathscr{L}_F(f)(x^q) \equiv 0 \mod \pi$$

\noindent
which is only possible if $\mathscr{L}_F(f) \in \pi\mathcal{O}_K[[x]]$. This completes the proof of the lemma.
\\
\\

\begin{lemma}
    If $\alpha_i \in \pi^{n-i}\textfrak{p}_0\mathcal{O}_{K_i}$ for $0 \le i \le n < \infty$, then there exists an $f\in \mathcal{O}_K[[x]]$ such that $f(u_i) = \alpha_i$ and $f(0) = 0$. Here $\textfrak{p}_0$ denotes the maximal ideal of $\mathcal{O}_{K_0}$. This is lemma 9 from \cite{Coleman1}.
\end{lemma}

\noindent
Proof: this follows from the observation that if 

$$g_{n,k} = \frac{ [\pi^{n+1}]\cdot [\pi^k]}{[\pi^{k+1}]}$$

\noindent
for $0 \le k \le n$, then $g_{n,k} \in \mathcal{O}_K[[x]]$ and

$$g_{n,k}(u_i) = 0$$

\noindent
for $0 \le i \le n$ and $i \neq k$. Also $g_{n,k}(u_k) = \pi^{n-k}u_0$. 
\\
\\
To get that $g_{n,k}(x) \in \mathcal{O}_K[[x]]$ write

$$\frac{[\pi^{n+1}]\cdot [\pi^k]}{[\pi^{k+1}]} = \frac{ [\pi^{n-k}]([\pi^{k+1}](x)) \cdot [\pi^k](x)}{[\pi^{k+1}](x)}$$

\noindent
Then note $x \mid [\pi^{n-k}](x)$ in $\mathcal{O}_K[[x]]$ which implies $g_{n,k}(x) \in \mathcal{O}_K[[x]]$. 
\\
\\
Next we need to show $g_{n,k}(u_i) = 0$ for $0 \le i \le n$ and $i \neq k$. To see this note that if $i < k$ then we have $[\pi^k](u_i) = 0$ and 

$$\frac{[\pi^{n+1}](x)}{[\pi^{k+1}](x)} \in \mathcal{O}_K[[x]]$$

\noindent
implying $g_{n,k}(u_i) = 0$ for such $i$. If $n \ge i > k$ then we can write

$$\frac{[\pi^{n+1}](x)}{[\pi^{k+1}](x)} = \frac{[\pi^{n-k}]([\pi^{k+1}](x))}{[\pi^{k+1}(x)]}$$

\noindent
Because $x \mid [\pi^{n-k}](x)$ in $\mathcal{O}_K[[x]]$, we have that the above is some series in $\mathcal{O}_K[[x]]$. Evaluating at $u_i$ when $n \ge i > k$ gives

$$\frac{[\pi^{n-k}](u_{i-k-1})}{u_{i-k-1}} = 0$$

\noindent
because $u_{i-k-1} \neq 0$. This implies that $g_{n,k}(u_i) = 0$ for such $i$. Note also that

$$g_{n,k}(u_k) = \pi^{n-k}u_0$$

\noindent
To see the above note that $[\pi^k](u_k) = u_0$. Note also that

$$\frac{[\pi^{n+1}](x)}{[\pi^{k+1}](x)} = \frac{[\pi^{n-k}]([\pi^{k+1}](x))}{[\pi^{k+1}(x)]}$$

\noindent
and $[\pi^{k+1}](u_k) = 0$, so that evaluating 

$$\frac{[\pi^{n+1}](x)}{[\pi^{k+1}](x)}$$

\noindent
at $u_k$ gives the constant term of the series 

$$\frac{[\pi^{n-k}](x)}{x}$$

\noindent
One can check that the constant term of this series is $\pi^{n-k}$.
\\
\\
Now we are given that each $\alpha_i \in \pi^{n-i}\textfrak{p}_0\mathcal{O}_{K_i}$, implying that for each $i$ there exists a series $f_i(x) \in \mathcal{O}_K[[x]]$ such that 

$$f_i(u_i) = \frac{\alpha_i}{\pi^{n-i}u_0}$$

\noindent
We get the above because

$$\frac{\alpha_i}{\pi^{n-i}u_0} \in \mathcal{O}_{K_i}$$

\noindent
and $u_i$ is a uniformizer for $\mathcal{O}_{K_i}$. We also use here that there exists a system of representatives $\{a_i\}$ for $\mathcal{O}_{K_i}/u_i\mathcal{O}_{K_i}$ with each $a_i \in \mathcal{O}_K$ because $K_i$ is totally ramified over $K$. One can then check that

$$f(x) = \sum_{i=0}^n f_i(x)g_{n,i}(x)$$

\noindent
satisfies $f(u_i) = \alpha_i$ for each $i$ with $0 \le i \le n$. Note also that $g_{n,i}(0) = 0$ for each $i$ with $0 \le i \le n$, so we get that $f(0) = 0$.

\pagebreak

\end{document}